\documentclass[a4paper,12pt]{amsart}
\usepackage{amssymb}
\usepackage{ifthen}
\usepackage[dvips]{graphicx}
\usepackage{cite}
\nonstopmode \numberwithin{equation}{section}
\usepackage{amssymb}
\usepackage{ifthen}
\usepackage{graphicx}
\usepackage{amsmath}
\usepackage[T1]{fontenc} 
\usepackage[utf8]{inputenc}
\usepackage[usenames,dvipsnames]{color}
\usepackage{color}
\usepackage[english]{babel}
\usepackage{fancyhdr}
\usepackage{fancybox}
\usepackage{tikz}

\nonstopmode \numberwithin{equation}{section}
\theoremstyle{plain}

\newtheorem{prop}{Proposition}

\newtheorem{conj}{Conjecture}

\theoremstyle{definition}
\newtheorem{defn}{Definition}[section]

\newtheorem{thm}{Theorem}[section]
\newtheorem{cor}{Corollary}[section]
\newtheorem{lem}{Lemma}[section]
\newtheorem{prob}{Problem}
\newtheorem{rem}{Remark}[section]
\newtheorem{ques}{Question}[section]

\newcounter{minutes}
\divide\time by 60
\newcounter{hours}
\multiply\time by 60
\addtocounter{minutes}{-\time}

\newcounter {own}
\def\theown {\thesection       .\arabic{own}}

\newenvironment{pf}[1][]{%
 \vskip 3mm
 \noindent
 \ifthenelse{\equal{#1}{}}%
  {{\slshape Proof. }}%
  {{\slshape #1.} }%
 }%
{\qed\bigskip}

\newcounter{alphabet}

\def\be{\begin{equation}}
\def\ee{\end{equation}}

\newcommand{\bee}{\begin{enumerate}}
\newcommand{\eee}{\end{enumerate}}

\newcommand{\blem}{\begin{lem}}
\newcommand{\elem}{\end{lem}}
\newcommand{\bthm}{\begin{thm}}
\newcommand{\ethm}{\end{thm}}
\newcommand{\bcor}{\begin{cor}}
\newcommand{\ecor}{\end{cor}}
\newcommand{\beg}{\begin{examp}}
\newcommand{\eeg}{\end{examp}}
\newcommand{\begs}{\begin{examples}}
\newcommand{\eegs}{\end{examples}}

\newcommand{\bdefn}{\begin{defn}}
\newcommand{\edefn}{\end{defn}}

\newcommand{\bprob}{\begin{prob}}
\newcommand{\eprob}{\end{prob}}
\newcommand{\bei}{\begin{itemize}}
\newcommand{\eei}{\end{itemize}}

\newcommand{\bcon}{\begin{conj}}
\newcommand{\econ}{\end{conj}}
\newcommand{\bcons}{\begin{conjs}}
\newcommand{\econs}{\end{conjs}}
\newcommand{\bprop}{\begin{prop}}
\newcommand{\eprop}{\end{prop}}
\newcommand{\br}{\begin{rem}}
\newcommand{\er}{\end{rem}}
\newcommand{\brs}{\begin{rems}}
\newcommand{\ers}{\end{rems}}
\newcommand{\bo}{\begin{obser}}
\newcommand{\eo}{\end{obser}}
\newcommand{\bos}{\begin{obsers}}
\newcommand{\eos}{\end{obsers}}
\newcommand{\bpf}{\begin{pf}}
\newcommand{\epf}{\end{pf}}
\newcommand{\ba}{\begin{array}}
\newcommand{\ea}{\end{array}}
\newcommand{\beq}{\begin{eqnarray}}
\newcommand{\beqq}{\begin{eqnarray*}}
\newcommand{\eeq}{\end{eqnarray}}
\newcommand{\eeqq}{\end{eqnarray*}}

\begin{document}

\title{Transcendental entire solutions of several general quadratic type PDEs and PDDEs in $ \mathbb{C}^2 $}

\author{Molla Basir Ahamed}
\address{Molla Basir Ahamed,
	Department of Mathematics,
	Jadavpur University,
	Kolkata-700032, West Bengal, India.}
\email{mbahamed.math@jadavpuruniversitry.in}

\author{Sanju Mandal}
\address{Sanju Mandal,
Department of Mathematics,
Jadavpur University,
Kolkata-700032, West Bengal, India.}
\email{sanjum.math.rs@jadavpuruniversity.in}

\subjclass[{AMS} Subject Classification:]{Primary 39A45, 30D35, 35M30, 32W50}

\def\thefootnote{}
\footnotetext{ {\tiny File:~\jobname.tex,
printed: \number\year-\number\month-\number\day,
          \thehours.\ifnum\theminutes<10{0}\fi\theminutes }
} \makeatletter\def\thefootnote{\@arabic\c@footnote}\makeatother

\begin{abstract}
The functional equations $ f^2+g^2=1 $ and $ f^2+2\alpha fg+g^2=1 $ are respectively called Fermat-type binomial and trinomial equations. It is of interest to know about the existence and form of the solutions of general quadratic functional equations. Utilizing Nevanlinna's theory for several complex variables, in this paper, we study the existence and form of the solutions to the general quadratic partial differential or partial differential-difference equations of the form $ af^2+2\alpha fg+b g^2+2\beta f+2\gamma g+C=0 $ in $ \mathbb{C}^2 $. Consequently, we obtain certain corollaries of the main results of this paper concerning binomial equations which generalize many results in [\textit{Rocky Mountain J. Math.} \textbf{51}(6) (2021),  2217-2235] in the sense of arbitrary coefficients.
\end{abstract}

\maketitle
\pagestyle{myheadings}
\markboth{Molla Basir Ahamed and Sanju Mandal}{Transcendental entire solutions of several general quadratic type PDEs and PDDEs in $ \mathbb{C}^2 $}
\section{Introduction}
The study of the Fermat type functional equation has been the interesting subject in the field of complex analysis in connection with extensions of Nevanlinna theory. Finding the precise solutions of different Fermat-type equations in several complex variables now-a-days becomes a central objectives in the value distribution theory of Nevanlinna. As we know that, partial differential equations are occurring in various areas of applied mathematics, such as fluid mechanics, nonlinear acoustic, gas dynamics and traffic flow (see \cite{Courant-Hilbert-I-1962},\cite{Garabedian-W-1964}). However, in general, it is difficult to find entire and meromorphic solutions for a nonlinear partial differential equations (PDEs) and differential-difference equations (DDEs). The study of complex DDEs can be traced back to Naftalevich’s research \cite{Naftalevich-1976,Naftalevich-1980}. Naftalevich have investigated the
meromorphic solutions on complex DDEs with one complex variable by employing the operator theory and iteration method. Moreover, by employing Nevanlinna theory and the method of complex analysis, there were a number of literature focusing on the solutions of the PDEs and DDEs and many variants of them.\vspace{1.2mm}

In a number of articles, it is dealt with the finite order solutions to the Fermat-type binomial and trinomial equations in $ \mathbb{C} $ and $ \mathbb{C}^2 $. Moreover, it appears recently in the literature that the transcendental entire solutions of the Fermat-type $ f^2+g^2=1 $ (binomial type) or  $ f^2+2\alpha fg+g^2=1 $ (trinomial type) equations os system of equations in several complex variables have been studied (see \cite{Xu-RMJM-2021,Xu-DM-2022,Xu-Liu-Li-JMAA-2020} and references therein). However, no attempt has been made here to develop methodologies on finding solutions of the general quadratic functional equations of the form $ af^2+2\alpha fg+bg^2+2\beta f+2\gamma g+C=0 $ (a general setting of binomial and trinomial type equations) in $ \mathbb{C} $ and $ \mathbb{C}^2 $. In this paper, our main objective is to find solutions to the system of general quadratic equations in $ \mathbb{C}^2 $. In this paper, we develop methods using some transformations and discuss how to solve reducing this general quadratic equations into a binomial form. The estimate we obtain in the course of proof seems to be independent interest. In fact, our main results provide a natural characterization of the solutions of the equations of our investigation. \vspace{1.2mm}

\par It has always been a well-known and interesting problem to investigate the existence and form of solution for Fermat-type equations
\begin{align}\label{eq-1.1}
	f^m(z)+g^m(z)=1
\end{align}
regard as the Fermat diophantine equation $ x^m+y^m=1 $ over functional fields, where $ n\geq 2 $ is an integer. In $ 1995 $, A. Wiles \cite{Wiles & Ann. Math. 1995}, A. Wiles and R. Taylor \cite{Tailor & Wiles & 1995} pointed out that the Fermat equation $ x^m+y^m=1 $ does not
admit nontrivial solutions in rational numbers for $ m\geq 3 $, and does admit nontrivial rational solutions for $ m=2 $. In fact, more earlier, Montel \cite{Montel & Paris & 1927}, Iyer \cite{Iyer & J. Indian. Math. Soc. & 1939} and Gross \cite{Gross & Bull. Amer. & 1966} had established results discussing the existence of solutions for Fermat type functional equation \eqref{eq-1.1} and pointed out the following:
\begin{enumerate}
	\item[(i)] for $ m=2 $, the entire solutions of \eqref{eq-1.1} are $ f(z)=\cos(\varphi(z)) $ and $ g(z)=\sin(\varphi(z)) $, where $ \varphi $ is an entire function;
	\item[(ii)] for $ m>2 $, there are no non-constant entire solutions of \eqref{eq-1.1};
	\item[(iii)] for $ m=2 $, the meromorphic solutions of \eqref{eq-1.1} are of the form
	\begin{align*}
		f(z)=\frac{2\varphi(z)}{1+\varphi^2(z)}\;\; \mbox{and}\;\; g(z)=\frac{1-\varphi^2(z)}{1+\varphi^2(z)},
	\end{align*}
	where $ \varphi $ is a meromorphic function;
	\item[(iv)] for $ m=3 $, the meromorphic solutions are of the form
	\begin{align*}
		f(z)=\frac{1}{2\wp(h)}\left(1+\frac{\wp^{\prime}(h)}{\sqrt{3}}\right)\;\; \mbox{and}\;\; g(z)=\frac{\omega}{2\wp(h)}\left(1-\frac{\wp^{\prime}(h)}{\sqrt{3}}\right),
	\end{align*}
	where $ \omega^3=1 $ and $ \wp $ satisfies $ \left(\wp^{\prime}\right)^2=4\wp^3-1 $;
	\item[(v)] for $ m>3 $, there are no non-constant meromorphic solutions.
\end{enumerate}
In $ 1933 $, Cartan \cite{Carten-1933} first considered the Fermat Diophantine functional equations $ f^m+g^n=1 $,
and observed that all the entire solutions $ f$ and $ g $ must be constants in case of when $ mn > m+n $ ; a similar observations are made also Gundersen and Hayman in \cite{Gundersen-Heyman-2004-BLMS}. Since then, researchers have paid their considerable attention to this study and investigated the existence of entire or meromorphic solutions for Fermat-type functional equation of the form $ f^n+g^m=1$, where $f, g$ are in general, meromorphic functions and $m,n\in\mathbb{N}$ (see e.g., \cite{Gross & Bull. Amer. & 1966,Montel & Paris & 1927,Tang & Liao & 2007,Yang & 1970,Yang & Li & 2004}). Continuing the study, in $ 2004 $, Yang and Li \cite{Yang & Li & 2004} investigated to find the form of solutions to the equation \eqref{eq-1.1} in view of replacing $g(z)$ by $f^{\prime}(z)$ when $ m=2 $, and showed that the transcendental entire solutions of $f(z)^{2}+f^{\prime}(z)^2=1$ must of the form $ 	f(z)=\frac{1}{2}\left(Ae^{\alpha z}+\frac{1}{A}e^{-\alpha z}\right), $  where $A, \alpha$ are non-zero complex constants. In \cite{Yang & Li & 2004}, the authors also studied the existence of meromorphic solutions to the equation $ 	 f^{2}(z)+\left(a_nf^{(n)}(z)+a_{n+1}f^{(n+1)}(z)\right)^2=1 $ for $a_n$ and $a_{n+1}$ be non-zero constants and proved that it has no transcendental meromorphic solutions. In $ 2015 $, Liu and Dong \cite{Liu & Dong & EJDE & 2015} further studied the existence of solutions of the following Fermat-type differential equation \begin{align}\label{eq-1.2}
	 (f(z)+f^{\prime}(z))^2+(f(z)+f^{\prime\prime}(z))^2=1
\end{align}
 and proved that equation \eqref{eq-1.2} does not admit any transcendental meromorphic solution. For different aspects of the solutions of Fermat-type functional equation, we refer to the articles \cite{Chen-JMAA-2022,Gundersen-CMFT-2017,Gundersen-PEMS-2020,Lu-Li & JMAA & 2019} and references therein. Moreover, for other aspects of Fermat-type functional equations and several questions on it, we refer the article \cite{Gundersen-CMFT-2017} by Gundersen \vspace{1.2mm}
 
For analogue study of Fermat-type functional equations in several complex variables, we recall here the classic result which plays a significant role.
\begin{thm}\cite{Saleeby-AM-2013}\label{th-1.1}
	For $ h:\mathbb{C}^n\rightarrow \mathbb{C} $ entire, the solutions of the equation $ f^2 +g^2=1 $ are characterize as follows:
	\begin{enumerate}
		\item [(i)] the entire solutions are $ f=\cos(h),\;g=\sin(h)$;
		\item[(ii)] the meromorphic solutions are of the form $ f=\frac{1-\beta^2}{1+\beta^2},\;g=\frac{2\beta}{1+\beta^2} $, with $ \beta $ being meromorphic on $ \mathbb{C}^n $.
	\end{enumerate}
\end{thm}
The first order non-linear partial differential equation $ u_{z_1}^2+u_{z_2}^2=1 $ is a special form of the functional equation $ f^2+g^2=1 $ in $ \mathbb{C}^2 $ that admits non-constant complex analytic solutions such as
\begin{align*}
	\begin{cases}
		f(z_1, z_2)=\cos(h(z_1, z_2))\\
		g(z_1, z_2)=\sin(h(z_1, z_2))
	\end{cases}
	\mbox{and} \;\;
	\begin{cases}
		f(z_1, z_2)=\sec(h(z_1, z_2))\\
		g(z_1, z_2)=i\tan(h(z_1, z_2))
	\end{cases}
\end{align*}
for an entire function $ h $ in $ \mathbb{C}^2 $.\vspace{1.2mm}

\noindent From now on, we will understand by $z+w=(z_1+w_1,z_2+w_2)$ for any $z=(z_1,z_2)$, $w=(w_1,w_2)$ are in $\mathbb{C}^2$, the shift of $f(z)$ is defined by $f(z+c)$, whereas the difference of $f(z)$ is defined by $\Delta_cf(z)=f(z+c)-f(z)$ (see \cite{Korhonen & CMFT & 2012}).\vspace{1.2mm}

Let us also recall here the definitions of difference equations and partial differential-difference equations in $ \mathbb{C}^n $.
\begin{defn}\cite{Zheng-Xu & Analysis math & 2021}
	An equation is called a difference equation, if the equation includes shifts or difference of $ f(z) $ in in $ \mathbb{C}^n $, which can be called $ DE $ for short. An equation is called a complex partial differential-difference equation, if this equation includes partial derivatives, shift or difference of $ f(z) $ in $ \mathbb{C}^n $, which can be called $ PDDE $ for short.
\end{defn}
Due to the development of the difference analogues of Nevanlinna theory in $ \mathbb{C} $ by Halburd and Korhonen \cite{Halburd & Korhonen & 2006}, Chiang and Feng \cite{Chiang & Feng & 2008} as well as in $ \mathbb{C}^n $  \cite{Cao-MN-2013,Cao-Korhonen-JMAA-2016,Gauthier-AM-1972, Vitter-DMJ-1977, Ye-MZ-1996}, many scholars paid their considerable attention to exploring several properties of the solutions to some complex difference equations both in $ \mathbb{C} $ and $ \mathbb{C}^n $. Most noticeably, $ 1995 $, Khavinson \cite{Khavinson & Am. Math. Mon & 1995} (also, Saleby \cite{Saleeby-A-1999}) derived that any entire solution of the partial differential equation $$\left(\frac{\partial u}{\partial z_1}\right)^2+\left(\frac{\partial u}{\partial z_2}\right)^2=1$$ must be linear, i.e., the solution takes the form $u(z_1, z_2)=az_1+bz_2+c$, where $a,b,c\in\mathbb{C}$, and $a^2+b^2=1$. In fact, Khavinson’s work inspired researchers on finding the entire solutions to the equations $ u_{z_1}^2+u_{z_2}^2=1 $ in $ \mathbb{C}^2 $. For example, Li \cite{Li & Nagoya & 2005,Li-TAMS-2005,Li & Arch. Math. & 2007} investigated on the partial differential equations with more general forms such as $u_{z_1}^2+u_{z_2}^2=p$, $u_{z_1}^2+u_{z_2}^2=e^q$, etc, where $p, q$ are polynomials in $\mathbb{C}^2$. With the help of the difference Nevanlinna theory for several complex variables (see \cite{Cao-Korhonen-JMAA-2016,Korhonen & CMFT & 2012}) obtained some interesting results on the existence of the solutions for some Fermat-type partial differential-difference equations from one complex variable to several several complex variables. \vspace{1.5mm}

In the recent years, Xu and his coauthors made some vital contribution establishing results on the Fermat-type equations with various settings in two complex variables. For extensive research work of Xu, the readers are referred to the articles \cite{Xu-RMJM-2021,Xu-AMP-2022,Xu-DM-2022,Xu-Liu-Li-JMAA-2020,Xu-Tu-Wang-RM-2021} and references therein. Moreover, for entire solutions of several quadratic binomial and trinomial partial differential-difference equations in $ \mathbb{C}^2 $, we refer to the article \cite{Haldar-Ahamed-AAMP-2022} and references there. In addition, entire solutions of certain type of nonlinear differential-difference equations are discussed in \cite{Chen-Nguyen-RMJ-2022}. Our interests in this paper to explore the solutions of equations considered in \cite{Xu-Tu-Wang-RM-2021} with a more general settings, namely for general quadratic functional equations. Henceforth, we first recall here a series of results obtained by Xu \emph{et al.} \cite{Xu-Tu-Wang-RM-2021} who have investigated the form of entire solutions to quadratic binomial functional equations and obtained the following result.

\begin{thm}\cite{Xu-Tu-Wang-RM-2021}\label{th-1.2}
Let $ f(z_1,z_2) $ be a transcendental entire solution with finite order of the partial differential equations $ \left(f(z) + \frac{\partial f(z)}{\partial z_1}\right)^2 + \left(f(z) + \frac{\partial f(z)}{\partial z_2}\right)^2 = 1. $ Then $ f(z_1,z_2) $ is of the form $ f(z_1,z_2)=\pm\frac{\sqrt{2}}{2} + \eta e^{-(z_1 + z_2)} $
or $ 	f(z_1,z_2)= \frac{1}{2}\sin(z_2 -z_1 + \eta_1) + \frac{1}{2}\cos(z_2 -z_1 + \eta_1) + \eta_2 e^{-(z_1 + z_2)} $, where $ \eta, \eta_1, \eta_2 \in\mathbb{C} $.
\end{thm}   
\begin{thm}\cite{Xu-Tu-Wang-RM-2021}\label{th-1.3}
The PDE	$ \left(f(z) + \frac{\partial f(z)}{\partial z_1}\right)^2 + \left(f(z) + \frac{\partial^2 f(z)}{\partial z^2_1}\right)^2 = 1. $
does not admit any transcendental entire solution with finite order.	
\end{thm}
\begin{thm}\cite{Xu-Tu-Wang-RM-2021}\label{th-1.4}
Let $ f(z_1,z_2) $ be a transcendental entire solution with finite order of the PDE 	$ \left(f(z) + \frac{\partial f(z)}{\partial z_1}\right)^2 + \left(f(z) + \frac{\partial^2 f(z)}{\partial z_1\partial z_2}\right)^2 = 1. $ Then $ f(z_1,z_2) $ is of the form $ f(z_1,z_2)=\pm\frac{\sqrt{2}}{2} + \eta e^{z_2 -z_1}, $ where $ \eta\in\mathbb{C} $.
\end{thm}

\begin{thm}\cite{Xu-Tu-Wang-RM-2021}\label{th-1.5}
Let $ c=(c_1,c_2)(\neq(0,0))\in\mathbb{C}^2 $ and $ s_0 = c_1 + c_2 $, and let $ f(z_1,z_2) $ be a transcendental entire solution with finite order of the PDE $ \left(f(z+c) + \frac{\partial f(z)}{\partial z_1}\right)^2 + \left(f(z+c) + \frac{\partial f(z)}{\partial z_2}\right)^2 = 1. $
Then $ f(z_1,z_2) $ is of the form
\begin{enumerate}
	\item [(i)] $ f(z_1,z_2)=\phi(z_1 + z_2), $
    where $ \phi(s) $ is a transcendental entire function with finite order in $ s:=z_1 + z_2 $ satisfying $ \phi(s + s_0) + \phi^{\prime}(s) =\pm\frac{\sqrt{2}}{2}; $
    \item[(ii)] $ f(z_1,z_2)=\frac{1 + i}{2(\alpha_1 -\alpha_2)}e^{L(z) + B} - \frac{1 -i}{2(\alpha_1-\alpha_2)}e^{-L(z)-B} + \phi(z_1+z_2), $
    where $ L(z)=\alpha_1z_1 + \alpha_2z_2 $, $ \alpha_1,\alpha_2,B\in\mathbb{C} $ and $ \phi(s) $ satisfies $ 	\alpha^2_1 + \alpha^2_2= -2,\; e^{L(c)}=-\frac{i\alpha_1 + \alpha_2}{1+i}=-\frac{1-i}{i\alpha_1-\alpha_2},\; \phi^{\prime}(s) + \phi(s+s_0)=0. $
\end{enumerate}
\end{thm}
\begin{thm}\cite{Xu-Tu-Wang-RM-2021}\label{th-1.6}
Let $ c=(c_1,c_2)(\neq(0,0))\in\mathbb{C}^2 $. Suppose that $ f(z_1,z_2) $ be a transcendental entire solution with finite order of the PDE $ \left(f(z+c) + \frac{\partial f(z)}{\partial z_1}\right)^2 + \left(f(z+c) + \frac{\partial^2 f(z)}{\partial z^2_1}\right)^2 = 1. $ Then $ f(z_1,z_2) $ is of the form	$ 	f(z_1,z_2)=\pm\frac{\sqrt{2}}{2} + \eta e^{z_1 +\frac{2k\pi i \pm \pi i -c_1}{c_2}z_2} $, where $ k\in\mathbb{Z}_{+} $ and $ \eta\in\mathbb{C} $; or $ f(z_1,z_2)=-\frac{1 + i}{2\alpha_1(\alpha_1 -1)}e^{L(z) + B} - \frac{1 -i}{2\alpha_1(\alpha_1-1)}e^{-L(z)-B} + \eta e^{z_1 +\frac{2k\pi i \pm \pi i -c_1}{c_2}z_2}, $ where $ L(z)=\alpha_1z_1 + \alpha_2z_2 $, and $ \eta\in\mathbb{C} $, $ \alpha_1,\alpha_2,c_1,c_2 $ satisfies $ (i\alpha_1 + \alpha^2_1)^2 =2,\; e^{2L(c)} =-i. $
\end{thm}
We see that all the equations considered in Theorems \ref{th-1.2} to \ref{th-1.6} are of the form $ F^2+G^2=1 $ (binomial type) and it is easy to see that $ F^2+2\alpha FG+G^2=1 $ (trinomial type)  is a  general setting of binomial, and moreover, the general quadratic equation of the form $ aF^2+2\alpha FG+bG^2+2\beta F+2\gamma G+C=0 $ is the natural extensions of both the binomial and trinomial equations. To the best of knowledge of the auhtors, since no investigation is done on exploring the existence and form of the solutions with this general quadratic functional equations. Thus, for further investigations, it is natural to raise the following question.
\begin{ques}\label{q-1.1}
Can we derive results corresponding to the equations in Theorems \ref{th-1.2} to \ref{th-1.6} considering general quadratic equations? If so, then can we find the precise form of the solutions of such equations?
\end{ques}
In this paper, our main aim is to answer Question \ref{q-1.1} completely. Before start our investigations, we want to mention here some results of Saleeby obtained for trinomial quadratic equations. In fact, Saleeby \cite{Saleeby-AM-2013} made some elementary observation on right factors of meromorphic function to describe complex analytic solutions to the quadratic functional equation $ f^2+2\alpha fg +g^2 =1 $, $ \alpha^2\neq1$ constant in $ \mathbb{C} $ and obtain a result associate with the partial differential equations $ u^2_x+2\alpha u_x u_y +u^2_y=1 $. Later in $ 2016 $, Liu and Yang \cite{Liu & Yang & 2016} studied the existence and the form of solutions of some quadratic trinomial functional equations and proved that if $ \alpha\neq \pm 1,0 $, then equation $ f(z)^2+2\alpha f(z)f^{\prime}(z)+f^{\prime}(z)^2=1 $ has no transcendental meromorphic solutions, whereas the equation $ f(z)^2+2\alpha f(z)f^{\prime}(z+c)+f^{\prime}(z+c)^2=1  $ must have transcendental entire solutions of order equal to one. In the recent years, trinomial $ PDDEs $ for functions of two complex variables are studied extensively. \vspace{1.2mm} 

Motivated by the techniques solving trinomial equations used by Saleeby in \cite{Saleeby-AM-2013}, our primary concern in this paper to establish results corresponding to Theorems \ref{th-1.2} to \ref{th-1.6} considering arbitrary coefficients in the equations. More precisely, we aim to generalize these results not only for trinomial cases but for complete general quadratic equations. In order to establish results corresponding to the equations in Theorems \ref{th-1.2} to \ref{th-1.6} for simplicity, we consider the general quadratic difference equations $ P(L_1(f), L_2(f))=0 $, $ P(L_1(f), L_3(f))=0 $ and  $ P(L_1(f), L_4(f))=0 $  or partial differential-difference equations $ P(M_1(f), M_2(f))=0 $ and $ P(M_1(f), M_3(f))=0 $  in $ \mathbb{C}^2 $.  Our main interest will be to find the solution of equation $ P(x, y)=0 $, where $ P(x, y)=ax^2+2\alpha xy+by^2+2\beta x+2\gamma y+C=0 $ which is irreducible \textit{i.e.,} equivalent to saying that the corresponding curve is non-singular. We define the following notations 
\begin{align*}
	\begin{cases}
		L_1(f(z))= f(z) + \dfrac{\partial f(z)}{\partial z_1}, L_2(f(z))= f(z) + \dfrac{\partial f(z)}{\partial z_2},\\ L_3(f(z)) = f(z) + \dfrac{\partial^2 f(z)}{\partial z^2_1},  L_4(f(z))= f(z) + \dfrac{\partial^2 f(z)}{\partial z_1 \partial z_2}, \\ M_1(f(z))= f(z+c) + \dfrac{\partial f(z)}{\partial z_1}, M_2(f(z)) = f(z+c) + \dfrac{\partial f(z)}{\partial z_2}, M_3(f(z)) = f(z+c) + \dfrac{\partial^2 f(z)}{\partial z^2_1}
	\end{cases}
\end{align*}
Employing elementary linear transformations and utilizing Nevanlinna's theory, in this paper, we will examine the existence and the precise form of the solutions of the following general quadratic $ PDDEs $.
\begin{align*}
\begin{cases}
aL_1(f)^2 + 2\alpha L_1(f)L_2(f)+ bL_2(f)^2 + 2\beta L_1(f) + 2\gamma L_2(f) + C  = 0, \\ aL_1(f)^2 + 2\alpha L_1(f)L_3(f)+ bL_3(f)^2 + 2\beta L_1(f) + 2\gamma L_3(f) + C  = 0, \\ aL_1(f)^2 + 2\alpha L_1(f)L_4(f)+ bL_4(f)^2 + 2\beta L_1(f) + 2\gamma L_4(f) + C  = 0, \\ aM_1(f)^2 + 2\alpha M_1(f)M_2(f)+ bM_2(f)^2 + 2\beta M_1(f) + 2\gamma M_2(f) + C  = 0, \\ aM_1(f)^2 + 2\alpha M_1(f)M_3(f)+ bM_3(f)^2 + 2\beta M_1(f) + 2\gamma M_3(f) + C  = 0.
\end{cases}
\end{align*}
\par In the above-mentioned equations, we assume the corresponding companion matrices are no-singular, i.e., $ \Delta\neq 0 $ , otherwise, each equation can be factorized linearly, which will eventually lead to the fact that the function $ f $ is a constant. Before we state the main results of this paper, below we discuss some techniques regarding how to transform a general quadratic equation into a binomial equation using transformations. 
\section{Finding solutions to the general quadratic equation with the help of some transformations}
The most general quadratic first order partial differential equation $ P(u_x, u_y)=0 $ over $ \mathbb{C} $ can be written as
\begin{equation}\label{eq-2.1}
	P(u_x, u_y):=au_x^2+2\alpha u_xu_y+bu_y^2+2\beta u_x+2\gamma u_y+C=0,
\end{equation}
where $ a, \alpha, b, \beta, \gamma, c $ are complex constants and $ (x, y)\in\mathbb{C}^2 $. It is easy to see that \eqref{eq-2.1} is associated with the algebraic equation
\begin{align}\label{eq-2.2}
	ax^2+2\alpha xy+by^2+2\beta x+2\gamma y+C=0.
\end{align}
This will represent a pair of straight lines only if the L.H.S of \eqref{eq-2.2}
can be resolve into two linear factors. Arranging the equation as a quadratic in $ x $, we have
\begin{align*}
	ax^2+2(\alpha y+\beta)x+(by^2+2\gamma y+C)=0.
\end{align*}
Therefore, if $ a\neq 0 $ a simple computation shows that
\begin{align*}
	x=\frac{-(\alpha y+\beta)\pm\sqrt{(\alpha y+\beta)^2-a(by^2+2\gamma y+C)}}{a}.
\end{align*}
Hence left hand side of \eqref{eq-2.2} can be resolve into two linear factors, if $ (\alpha y+\beta)^2-a(by^2+2\gamma y+C) $, i.e, $ y^2(\alpha^2 -ab)+2(\alpha\beta-a\gamma)y+(\beta^2 -aC) $ would be a perfect square. This is a quadratic in $ y $, and hence is a perfect square, if its discriminant be zero, i.e, if $ 4(\alpha\beta-a\gamma)^2 -4(\alpha^2 -ab)(\beta^2 -aC)=0 $\;\; or \;\; 
\begin{align}\label{eq-2.3}
	a(abC+2\alpha\beta\gamma-a\gamma^2-b\beta^2-C\alpha^2)=0.
\end{align} The companion matrix of the quadratic equation \eqref{eq-2.2} is the following
\[
A=\begin{bmatrix}
	a & \alpha & \beta\\
	\alpha & b & \gamma\\
	\beta & \gamma & C
\end{bmatrix}
.\]
Let $ \Delta:=\det(A) $, where $ \Delta $ is known to be an invariant of the conic as $ \Delta=0 $ in invariant under linear transformations (see \cite{Gibson-CUP_1998}). A simple computation shows that $ \Delta=abC+2\alpha\beta\gamma-a\gamma^2-b\beta^2-C\alpha^2 $. Now since $ a\neq0 $, from \eqref{eq-2.3}, we see that $\Delta=0 $. Therefore, when $ \Delta=0 $ the equation \eqref{eq-2.2} represent a pair of straight lines. These pair of straight lines will be parallel if $ \alpha^2=ab $ and be intersecting if $ \alpha^2\neq ab $. \vspace{1.2mm}
In addition, when $ \Delta\neq 0 $, the equation \eqref{eq-2.2} is transformed into a circle $ x^2+y^2=1 $, or, an ellipse $ x^2/a^2+y^2/b^2=1 $, or, a parabola, $ y^2=x $. Define, $ D=\det\begin{bmatrix}
	a & \alpha\\
	\alpha & b
\end{bmatrix}=ab-\alpha^2 $. Let $ \alpha^2\neq ab $. Then 
under the transformation 
\begin{align}\label{eq-2.4}
	x=X+x_1\; \mbox{and}\; y=Y+y_1, \; \mbox{where}\; x_1=\frac{\alpha\gamma-b\beta}{ab-\alpha^2}\; \mbox{and}\; y_1=\frac{\alpha\beta-a\gamma}{ab-\alpha^2}
\end{align}
a straightforward computation shows that equation $\eqref{eq-2.4}  $ reduces to
\begin{align}\label{eq-2.5}
	aX^2+2\alpha XY +bY^2+C_1=0,
\end{align}
where $ C_1=ax^2_1+2\alpha x_1 y_1+by^2_1+2\beta x_1+2\gamma y_1 +C. $ In view of the values of $ x_1 $ and $ y_1 $, we see that $ C_1 $ can be written as $ 	C_1=\frac{abC+2\alpha\beta\gamma-a\gamma^2-b\beta^2-C\alpha^2}{ab-\alpha^2}=\frac{\Delta}{D}. $ Hence, the equation $ \eqref{eq-2.5} $ can be written as
\begin{align}\label{eq-2.6}
	aX^2+2\alpha XY +bY^2+\frac{\Delta}{D}=0,
\end{align}
where $ D\neq 0 $ as $ \alpha^2\neq ab $. To remove the $ XY $-term from \eqref{eq-2.6}, we take the transformation $ (X, Y)=\left(\mathcal{U}\xi^{\pm}_1 - \mathcal{V}\eta^{\pm}_1,\; \mathcal{U}\eta^{\pm}_1+ \mathcal{V}\xi^{\pm}_1\right), $
where
\begin{align*}
	\begin{cases}
		\xi^{\pm}_1:=\dfrac{2\alpha}{\sqrt{\left((b-a)\pm\sqrt{(a-b)^2 +4\alpha^2}\right)^2 +4\alpha^2}},\vspace{1.2mm}\\ \eta^{\pm}_1:=\dfrac{(b-a)\pm\sqrt{(a-b)^2 +4\alpha^2}}{\sqrt{\left((b-a)\pm\sqrt{(a-b)^2 +4\alpha^2}\right)^2 +4\alpha^2}}.
	\end{cases}
\end{align*}
Under this transformation, the equation $ \eqref{eq-2.6} $ is transformed to
\begin{align}\label{eq-2.7}
	\left(\frac{\mathcal{U}}{1/\sqrt{A^{\pm}}}\right)^2+	\left(\frac{\mathcal{V}}{1/\sqrt{B^{\mp}}}\right)^2+\frac{\Delta}{D}=0,
\end{align}
where 
\begin{align*}
	A^{\pm}:=\dfrac{(a+b)\pm\sqrt{(a-b)^2 +4\alpha^2}}{2}\; \mbox{and}\; B^{\mp}:=\dfrac{(a+b)\mp\sqrt{(a-b)^2 +4\alpha^2}}{2}.
\end{align*}
It is easy to see that \eqref{eq-2.7} can be written as
\begin{align}\label{eq-2.8}
	\left(\sqrt{\frac{DA^{\pm}}{-\Delta}}\;\mathcal{U}\right)^2 + \left(\sqrt{\frac{DB^{\mp}}{-\Delta}}\;\mathcal{V}\right)^2 =1,
\end{align}
which is associated with a circle in $ \mathbb{C}^2 $. Therefore, in view of part -(i) of Theorem \ref{th-1.1}, we must have
\begin{align*}
	(\mathcal{U}, \mathcal{V})=\left(\frac{\cos h(z)}{\sqrt{\frac{DA^{\pm}}{-\Delta}}}, \frac{\sin h(z)}{\sqrt{\frac{DB^{\mp}}{-\Delta}}}\right),
\end{align*}
where $ h :\mathbb{C}^2\rightarrow \mathbb{C} $ is an entire function. In view of the above observations, it follows from $ \eqref{eq-2.4} $ that the transcendental solutions $ (x, y) $ to the general quadratic equation \eqref{eq-2.2} take the form
\begin{align}\label{eq-2.9}
	(x, y)=\left(D_{11}\cos h(z) -D_{12}\sin h(z) + T_1, E_{11}\cos h(z) +E_{12}\sin h(z) + T_2\right)
\end{align}
where
\begin{align}\label{eq-2.10}
	\begin{cases}
		D_{11}=\dfrac{\xi^{\pm}_1}{\sqrt{\frac{DA^{\pm}}{-\Delta}}},\;\; D_{12}=\dfrac{\eta^{\pm}_1}{\sqrt{\frac{DB^{\mp}}{-\Delta}}}\;\;\mbox{and}\;\; T_1=\dfrac{\alpha\gamma-b\beta}{ab-\alpha^2},\vspace{1.5mm}\\ E_{11}=\dfrac{\eta^{\pm}_1}{\sqrt{\frac{DA^{\pm}}{-\Delta}}},\;\; E_{12}=\dfrac{\xi^{\pm}_1}{\sqrt{\frac{DB^{\mp}}{-\Delta}}}\;\;\mbox{and}\;\; T_2=\dfrac{\alpha\beta-a\gamma}{ab-\alpha^2}.
	\end{cases}
\end{align}
\section{Main results and their consequences}
 In this paper, we investigate the existence and the form of entire solutions for general quadratic PDEs and PDDES with two complex variables.  We obtain Theorems \ref{th-3.1} to \ref{th-3.5} corresponding to Theorems \ref{th-1.2} to \ref{th-1.6} respectively, and hence Question \ref{q-1.1} is answered completely. The first main result of the paper is the following which extends Theorem \ref{th-1.4} for general quadratic PDEs.
\begin{thm}\label{th-3.1}
Let $ a, b, C, \alpha, \beta, \gamma\in\mathbb{C} $ with $\alpha^2\neq ab $, $ \Delta\neq 0 $ and $ f(z_1,z_2) $ be a finite order transcendental entire solutions for the general quadratic partial differential equation 
\begin{align}\label{eq-3.1}
	aL_1(f)^2 + 2\alpha L_1(f)L_2(f)+ bL_2(f)^2 + 2\beta L_1(f) + 2\gamma L_2(f) + C  = 0.
\end{align}
Then $ f(z_1,z_2) $ is of the forms
\begin{enumerate}
	\item [(i)] $ f(z_1,z_2)= A_0 + R_4 e^{-(z_1 +z_2)}, $ where $ K_1, K_2, R_4 $ are constants with $ K^2_1 + K^2_2= 1 $,
    \begin{align*}
    	\begin{cases}
    		A_0 = \dfrac{K_1\xi^{\pm}_1}{\sqrt{\frac{DA^{\pm}}{-\Delta}}} - \dfrac{K_2\eta^{\pm}_1}{\sqrt{\frac{DB^{\mp}}{-\Delta}}} + \dfrac{\alpha\gamma- b\beta}{ab- \alpha^2},\; 	K_2 =\dfrac{-R_2R_3\pm R_1\sqrt{R^2_1 +R^2_2 -R^2_3}}{R^2_1 + R^2_2},\vspace{2mm}\\
    		R_1 =\dfrac{\xi^{\pm}_1-\eta^{\pm}_1}{\sqrt{\frac{DA^{\pm}}{-\Delta}}},\; R_2=\dfrac{\xi^{\pm}_1+ \eta^{\pm}_1}{\sqrt{\frac{DB^{\mp}}{-\Delta}}}\; \mbox{and}\; R_3=\dfrac{\beta(\alpha + b)- \gamma(\alpha + a)}{ab-\alpha^2}.
    	\end{cases}
    \end{align*}
    \item[(ii)] \begin{align*}
    	f(z_1,z_2)&=\left(\dfrac{E_{11}+E_{12}}{2}\right)\sin\left(\dfrac{R_{12}}{R_{11}}z_1 + z_2 + R_5\right) \\& + \left(\dfrac{E_{11}-E_{12}}{2}\right)\cos \left(\dfrac{R_{12}}{R_{11}}z_1 + z_2 + R_5\right) + R_6 e^{-(z_1 + z_2)} + T_2.
    \end{align*} where $ R_5 $ is a constant and $ T_1 = T_2 $.
\end{enumerate}
\end{thm}
Before, we state some consequences of the main result as the precise form of the solutions for trinomial and binomial equations, we introduce some notations below. Denote
\begin{align*}
\begin{cases}
K_{11}=\sqrt{(a+b)\pm\sqrt{(a-b)^2 + 4\alpha^2}},\vspace{1.5mm}\\ K_{12}=\sqrt{\left((b-a)\pm\sqrt{(a-b)^2+4\alpha^2}\right)^2 +4\alpha^2},\vspace{1.5mm}\\ K_{13}=\sqrt{(a+b)\mp\sqrt{(a-b)^2 + 4\alpha^2}}, K_{14}=(b-a)\pm\sqrt{(a-b)^2+4\alpha^2)},\\ A_{11}=\dfrac{2\sqrt{2}\alpha}{K_{11}K_{12}}, A_{12}=\dfrac{\sqrt{2}K_{14}}{K_{12}K_{13}}, B_{11}=\dfrac{\sqrt{2}K_{14}}{K_{11}K_{12}}, B_{12}=\dfrac{2\sqrt{2}\alpha}{K_{12}K_{13}}.
\end{cases}
\end{align*}
We obtain the following result from Theorem \ref{th-3.1} showing the precise form of the solutions in case of trinomial quadratic PDEs.
\begin{cor}\label{cor-3.1}
Let $ a, b, C, \alpha\in\mathbb{C} $ with $\alpha^2\neq ab $, $ \Delta\neq 0 $ and $ f(z_1,z_2) $ be a finite order transcendental entire solutions for the general quadratic partial differential equation 
\begin{align*}
	a\left(f(z) + \dfrac{\partial f(z)}{\partial z_1}\right)^2 + 2\alpha \left(f(z) + \dfrac{\partial f(z)}{\partial z_1}\right)\left( f(z) + \dfrac{\partial f(z)}{\partial z_2}\right)+ b\left( f(z) + \dfrac{\partial f(z)}{\partial z_2}\right)^2 = 1.
\end{align*}
Then $ f(z_1,z_2) $ is of the forms
\begin{enumerate}
	\item [(i)] $ 	f(z_1,z_2)= A_{31} + R_4 e^{-(z_1 +z_2)}, $	where $ K_1, K_2, R_4 $ are constants with $ K^2_1 + K^2_2= 1 $, $ A_{31} = A_{11}K_1 - A_{12}K_2 $, $ K_2 =\pm\dfrac{S_{11}}{\sqrt{S^2_{11} + S^2_{12}}} $, where $ S_{11} =\frac{\sqrt{2}(2\alpha - K_{14})}{K_{11}K_{12}}\; \mbox{and}\; S_{12} =\frac{\sqrt{2}(2\alpha + K_{14})}{K_{11}K_{12}}. $
	\item[(ii)] \begin{align*}
		f(z_1,z_2)&=\left(\dfrac{B_{11}+B_{12}}{2}\right)\sin\left(\dfrac{S_{14}}{S_{13}}z_1 + z_2 + R_5\right) \\& + \left(\dfrac{B_{11}-B_{12}}{2}\right)\cos \left(\dfrac{S_{14}}{S_{13}}z_1 + z_2 + R_5\right) + R_6 e^{-(z_1 + z_2)}.
	\end{align*} where $ R_5 $ is a constant, $ S_{13}= B_{11}(A_{11}-B_{11}) - B_{12}(A_{12}+B_{12})$ and $ S_{14}= A_{11}(A_{11}-B_{11}) + A_{12}(A_{12}+B_{12}) $.
\end{enumerate}
\end{cor}
We obtain the following corollary of Theorem \ref{th-3.1} and it gives the solutions in generalized form of the same in Theorem \ref{th-1.2}.
\begin{cor}\label{cor-3.2}
Let $ a, b\in\mathbb{C} $ with $ \Delta\neq 0 $ and $ f(z_1,z_2) $ be a finite order transcendental entire solutions for the general quadratic partial differential equation 
\begin{align*}
a\left(f(z) + \dfrac{\partial f(z)}{\partial z_1}\right)^2 + b\left( f(z) + \dfrac{\partial f(z)}{\partial z_2}\right)^2 = 1.
\end{align*} 
Then $ f(z_1,z_2) $ is of the forms
\begin{enumerate}
\item [(i)] $ f(z_1,z_2)= A_0 + R_4 e^{-(z_1 +z_2)}, $ where $ K_{1*}, K_{2*}, R_4 $ are constants with $ K^2_{1*} + K^2_{2*}= 1 $, $ A_0 = - A_{21}K_2 $ and $ K_{2*} =\mp\frac{K_{21}}{\sqrt{K^2_{21}+K^2_{22}}} $
\item[(ii)] 
\begin{align*}
f(z_1,z_2)&=\dfrac{A_{22}}{2}\sin\left(-\dfrac{A^2_{21}}{A^2_{22}}z_1 + z_2 + R_5\right)\\& + \dfrac{A_{22}}{2}\cos \left(-\dfrac{A^2_{21}}{A^2_{22}}z_1 + z_2 + R_5\right) + R_6 e^{-(z_1 + z_2)}, \mbox{where}\; R_5,  R_6\; \mbox{are constants}.
\end{align*} 
\end{enumerate}
\end{cor}

\begin{thm}\label{th-3.2}
Let $ a, b, C, \alpha, \beta, \gamma\in\mathbb{C} $ with $\alpha^2\neq ab $ and $ \Delta\neq 0 $. The partial differential equation
\begin{align}\label{eq-3.2}
	aL_1(f)^2 + 2\alpha L_1(f)L_3(f)+ bL_3(f)^2 + 2\beta L_1(f) + 2\gamma L_3(f) + C  = 0
\end{align} 
does not admit any transcendental entire solution with finite order.   
\end{thm}
\begin{rem}
	In view of the proof of Theorem \ref{th-3.2}, it can be easily shown that the trinomial quadratic PDE, $ a\left(f(z) + \frac{\partial f(z)}{\partial z_1}\right)^2 + 2\alpha \left(f(z) + \frac{\partial f(z)}{\partial z_1}\right)\left(f(z) + \frac{\partial^2 f(z)}{\partial z^2_1}\right)+ b\left(f(z) + \frac{\partial^2 f(z)}{\partial z^2_1}\right)^2= 1 $ and also the binomial quadratic PDE, $ a\left(f(z) + \frac{\partial f(z)}{\partial z_1}\right)^2 + b\left(f(z) + \frac{\partial^2 f(z)}{\partial z^2_1}\right)^2= 1 $ do not admit any transcendental entire solutions. Consequently, the Theorem \ref{th-1.3} is established with more general setting of the equation.
\end{rem}
\begin{thm}\label{th-3.3}
Let $ a, b, C, \alpha, \beta, \gamma\in\mathbb{C} $ with $\alpha^2\neq ab $, $ \Delta\neq 0 $ and $ f(z_1,z_2) $ be a finite order transcendental entire solutions for the general quadratic partial differential equation 
\begin{align}\label{eq-3.3}
	aL_1(f)^2 + 2\alpha L_1(f)L_4(f)+ bL_4(f)^2 + 2\beta L_1(f) + 2\gamma L_4(f) + C  = 0.
\end{align}
Then $ f(z_1,z_2) $ is of the forms $ f(z_1,z_2)= -A_0 + B_{11}e^{(z_1 - z_2)} $, where $ K_5, K_6, B_{11} $ are constant with $ K^2_5 + K^2_6= 1 $,
\begin{align*}
\begin{cases}
A_0 = \dfrac{K_5\xi^{\pm}_1}{\sqrt{\frac{DA^{\pm}}{-\Delta}}} - \dfrac{K_6\eta^{\pm}_1}{\sqrt{\frac{DB^{\mp}}{-\Delta}}} + \dfrac{\alpha\gamma- b\beta}{ab- \alpha^2},\; 	K_6 =\dfrac{-R_2R_3\pm R_1\sqrt{R^2_1 +R^2_2 -R^2_3}}{R^2_1 + R^2_2},\\
R_1 =\dfrac{\xi^{\pm}_1-\eta^{\pm}_1}{\sqrt{\frac{DA^{\pm}}{-\Delta}}},\; R_2=\dfrac{\xi^{\pm}_1+ \eta^{\pm}_1}{\sqrt{\frac{DB^{\mp}}{-\Delta}}}\; \mbox{and}\; R_3=\dfrac{\beta(\alpha + b)- \gamma(\alpha + a)}{ab-\alpha^2}.
\end{cases}
\end{align*}
\end{thm}
In case of trinomial quadratic PDEs, we obtain the following result as a consequence of Theorem \ref{th-3.3}.
\begin{cor}
 Let $ a, b, C, \alpha\in\mathbb{C} $ with $\alpha^2\neq ab $, $ \Delta\neq 0 $ and $ f(z_1,z_2) $ be a finite order transcendental entire solutions for the general quadratic trinomial PDE
\begin{align*}
a\left(f(z) + \frac{\partial f(z)}{\partial z_1}\right)^2 + 2\alpha \left(f(z) + \frac{\partial f(z)}{\partial z_1}\right)\left(f(z) + \dfrac{\partial^2 f(z)}{\partial z_1 \partial z_2}\right)+ b\left(f(z) + \dfrac{\partial^2 f(z)}{\partial z_1 \partial z_2}\right)^2 =1.
\end{align*}
Then $ f(z_1,z_2) $ is of the forms $ f(z_1,z_2)= -A_{31} + B_{11}e^{(z_1 - z_2)} $, where $ K_{5*}, K_{6*}, B_{11} $ are constant with $ K^2_{5*} + K^2_{6*}= 1 $, $ A_{31} = A_{11}K_{5*} - A_{12}K_{6*} $ and $ K_{6*} =\pm\frac{S_{11}}{\sqrt{S^2_{11} + S^2_{12}}}, $ where $ 	S_{11} =\frac{\sqrt{2}(2\alpha - K_{14})}{K_{11}K_{12}}\; \mbox{and}\; S_{12} =\frac{\sqrt{2}(2\alpha + K_{14})}{K_{11}K_{12}}. $
\end{cor}
We obtain the following corollary of Theorem \ref{th-3.3} which generalizes Theorem \ref{th-1.4}.
\begin{cor}
	Let $ a, b\in\mathbb{C} $ with $ \Delta\neq 0 $ and $ f(z_1,z_2) $ be a finite order transcendental entire solutions for the general quadratic binomial PDE 
	\begin{align*}
		a\left(f(z) + \frac{\partial f(z)}{\partial z_1}\right)^2 + b\left(f(z) + \dfrac{\partial^2 f(z)}{\partial z_1 \partial z_2}\right)^2 =1.
	\end{align*}
	Then $ f(z_1,z_2) $ is of the forms $ f(z_1,z_2)= -A_0 + B_{11}e^{(z_1 - z_2)} $
	where $ K_{11*}, K_{22*}, B_{11} $ are constant with $ K^2_{11*} + K^2_{22*}= 1 $, $ A_0 = - A_{21}K_{22*} $ and $ 	K_{22*} =\mp\frac{K_{21}}{\sqrt{K^2_{21}+K^2_{22}}}. $
\end{cor}
\begin{thm}\label{th-3.4}
Let $ a, b, C, \alpha, \beta, \gamma\in\mathbb{C} $ with $\alpha^2\neq ab $, $ \Delta\neq 0 $, $ c=(c_1,c_2)(\neq (0,0))\in\mathbb{C}^2 $, $ s_0=c_1+c_2 $ and $ f(z_1,z_2) $ be a finite order transcendental entire solutions for the general quadratic partial differential equation
\begin{align}\label{eq-3.4}
	aM_1(f)^2 + 2\alpha M_1(f)M_2(f)+ bM_2(f)^2 + 2\beta M_1(f) + 2\gamma M_2(f) + C  = 0.
\end{align}	
Then $ f(z_1,z_2) $ is one of the following forms
\begin{enumerate}
	\item [(i)] $ 	f(z_1,z_2)= \psi(z_1 +z_2) $, where $ \psi(s) $ is a transcendental entire function with finite order in $ s:=z_1+z_2 $, satisfying $ \psi^{\prime}(s) + \psi(s+ s_0) = A_0, $
     where $ K_7, K_8 $ are constants with $ K^2_7 + K^2_8= 1 $, where
     \[
     \begin{cases}
     A_0 = \dfrac{K_7\xi^{\pm}_1}{\sqrt{\frac{DA^{\pm}}{-\Delta}}} - \dfrac{K_8\eta^{\pm}_1}{\sqrt{\frac{DB^{\mp}}{-\Delta}}} + \dfrac{\alpha\gamma- b\beta}{ab- \alpha^2},\; K_6 =\dfrac{-R_2R_3\pm R_1\sqrt{R^2_1 +R^2_2 -R^2_3}}{R^2_1 + R^2_2},\vspace{2mm}\\ R_1 =\dfrac{\xi^{\pm}_1-\eta^{\pm}_1}{\sqrt{\frac{DA^{\pm}}{-\Delta}}},\; R_2=\dfrac{\xi^{\pm}_1+ \eta^{\pm}_1}{\sqrt{\frac{DB^{\mp}}{-\Delta}}}\; \mbox{and}\; R_3=\dfrac{\beta(\alpha + b)- \gamma(\alpha + a)}{ab-\alpha^2}.
     \end{cases}
     \]
     \item[(ii)] \begin{align*}
     	f(z_1,z_2)=\dfrac{(D_{11} - E_{11})}{(\alpha_1 - \alpha_2)}\sin (L(z) + B) + \dfrac{(D_{12} + E_{12})}{(\alpha_1 - \alpha_2)}\cos (L(z)+ B) + (T_1 - T_2)z_1 + \psi_{1}(s).
     \end{align*}
     where 
     \begin{align*}
     	e^{2iL(c)}=\dfrac{(D_{11}+iD_{12})\alpha_2 - (E_{11}-iE_{12})\alpha_1}{(E_{11}+iE_{12})\alpha_1 - (D_{11}-iD_{12})\alpha_2}\; \mbox{and}\; L(z)=\alpha_1z_1 + \alpha_2z_2,\;  \alpha_1,\alpha_2, B\in\mathbb{C}
     \end{align*}
     and
     \begin{align*}
     	(E^2_{11} + E^2_{12})\alpha^2_1 &+ (D^2_{11} + D^2_{12})\alpha^2_2 - 2(D_{11}E_{11} - D_{12}E_{12})\alpha_1\alpha_2 \\&+ [i(D_{11}-E_{11})+(D_{12}+E_{12})]^2 = 0
     \end{align*}
      and $ \psi_{1}(s) $ satisfies 
      \begin{align*}
      	&\psi^{\prime}_{1}(s) + \psi_{1}(s+s_0)\\ &= \left(\dfrac{E_{11}\alpha_1 - D_{11}\alpha_2}{\alpha_1-\alpha_2}\right)\cos (L(z)+B) + \left(\dfrac{E_{12}\alpha_1 - D_{12}\alpha_2}{\alpha_1-\alpha_2}\right)\sin (L(z)+B)\\& -\left(\dfrac{D_{11}-E_{11}}{\alpha_1-\alpha_2}\right)\sin(L(z)+L(c)+B) - \left(\dfrac{D_{12}-E_{12}}{\alpha_1-\alpha_2}\right)\cos(L(z)+L(c)+B)\\& - (T_1-T_2)(z_1+c) + T_2.
      \end{align*}
\end{enumerate}
\end{thm} 
For trinomial quadratic PDEs, the following result is a corollary of Theorem \ref{th-3.4}.
\begin{cor}\label{cor-3.3}
Let $ a, b, C, \alpha\in\mathbb{C} $ with $\alpha^2\neq ab $, $ \Delta\neq 0 $, $ c=(c_1,c_2)(\neq (0,0))\in\mathbb{C}^2 $, $ s_0=c_1+c_2 $ and $ f(z_1,z_2) $ be a finite order transcendental entire solutions for the general quadratic partial differential equation
\begin{align*}
	a\left(f(z+c) + \frac{\partial f(z)}{\partial z_1}\right)^2 + 2\alpha \left(f(z+c) + \frac{\partial f(z)}{\partial z_1}\right)\left(f(z+c) + \frac{\partial f(z)}{\partial z_2}\right)+ b\left(f(z+c) + \frac{\partial f(z)}{\partial z_2}\right)^2 = 1.
\end{align*}	
Then $ f(z_1,z_2) $ is one of the following forms
\begin{enumerate}
	\item [(i)] $ 	f(z_1,z_2)= \psi(z_1 +z_2) $, where $ \psi(s) $ is a transcendental entire function with finite order in $ s:=z_1+z_2 $, satisfying $ \psi^{\prime}(s) + \psi(s+ s_0) = A_{31}, $ where $ K_{7*}, K_{8*} $ are constant with $ K^2_{7*} + K^2_{8*}= 1 $, $ A_{31} = A_{11}K_{7*} - A_{12}K_{8*} $, $ K_{8*} =\pm\frac{S_{11}}{\sqrt{S^2_{11} + S^2_{12}}} $, $ S_{11} =\frac{\sqrt{2}(2\alpha - K_{14})}{K_{11}K_{12}}\; \mbox{and}\; S_{12} =\frac{\sqrt{2}(2\alpha + K_{14})}{K_{11}K_{12}}. $
	\item[(ii)] \begin{align*}
		f(z_1,z_2)=\dfrac{(A_{11} - B_{11})}{(\alpha_1 - \alpha_2)}\sin (L(z) + B) + \dfrac{(A_{12} + B_{12})}{(\alpha_1 - \alpha_2)}\cos (L(z)+ B) + \psi_{1}(s).
	\end{align*}
	where
	\begin{align*}
		\begin{cases}
	L(z)=\alpha_1z_1 + \alpha_2z_2 $, $ \alpha_1,\alpha_2, B\in\mathbb{C},\;
	e^{2iL(c)}=\dfrac{(A_{11}+iA_{12})\alpha_2 - (B_{11}-iB_{12})\alpha_1}{(B_{11}+iB_{12})\alpha_1 - (A_{11}-iA_{12})\alpha_2},\vspace{2mm}\\
	(B^2_{11} + B^2_{12})\alpha^2_1 + (A^2_{11} + A^2_{12})\alpha^2_2 - 2(A_{11}B_{11} - A_{12}B_{12})\alpha_1\alpha_2 \\+ [i(A_{11}-B_{11})+(A_{12}+B_{12})]^2 = 0.
		\end{cases}
	\end{align*}
	and $ \psi_{1}(s) $ satisfies 
	\begin{align*}
		&\psi^{\prime}_{1}(s) + \psi_{1}(s+s_0)\\ &= \left(\dfrac{B_{11}\alpha_1 - A_{11}\alpha_2}{\alpha_1-\alpha_2}\right)\cos (L(z)+B) + \left(\dfrac{B_{12}\alpha_1 - A_{12}\alpha_2}{\alpha_1-\alpha_2}\right)\sin (L(z)+B)\\& -\left(\dfrac{A_{11}-B_{11}}{\alpha_1-\alpha_2}\right)\sin(L(z)+L(c)+B) - \left(\dfrac{A_{12}-B_{12}}{\alpha_1-\alpha_2}\right)\cos(L(z)+L(c)+B).
	\end{align*}
\end{enumerate}
\end{cor}
The following corollary of Theorem \ref{th-3.4} generalizes Theorem \ref{th-1.5}.
\begin{cor}\label{cor-3.6}
	Let $ a, b\in\mathbb{C} $ with $ \Delta\neq 0 $, $ c=(c_1,c_2)(\neq (0,0))\in\mathbb{C}^2 $, $ s_0=c_1+c_2 $ and $ f(z_1,z_2) $ be a finite order transcendental entire solutions for the general quadratic partial differential equation
	\begin{align*}
	a\left(f(z+c) + \frac{\partial f(z)}{\partial z_1}\right)^2 + b\left(f(z+c) + \frac{\partial f(z)}{\partial z_2}\right)^2 = 1.
	\end{align*}	
	Then $ f(z_1,z_2) $ is one of the following forms
	\begin{enumerate}
		\item [(i)] $ f(z_1,z_2)= \psi(z_1 +z_2) $,	where $ \psi(s) $ is a transcendental entire function with finite order in $ s:=z_1+z_2 $, satisfying $ \psi^{\prime}(s) + \psi(s+ s_0) = A_0, $
		where $ K_{7**},K_{8**} $ are constant with $ K^2_{7**} + K^2_{8**}= 1 $, $ A_0 = - A_{21}K_{8**} $ and $ K_{8**} =\mp\frac{K_{21}}{\sqrt{K^2_{21}+K^2_{22}}}. $
		\item[(ii)] 
		\begin{align*}
			f(z_1,z_2)=\dfrac{- A_{22}}{(\alpha_1 - \alpha_2)}\sin (L(z) + B) + \dfrac{A_{21}}{(\alpha_1 - \alpha_2)}\cos (L(z)+ B) + \psi_{1}(s).
		\end{align*}
		where $ L(z)=\alpha_1z_1 + \alpha_2z_2 $, $ \alpha_1,\alpha_2, B\in\mathbb{C} $ and $ e^{2iL(c)}=\dfrac{iA_{21}\alpha_2 - A_{22}\alpha_1}{iA_{21}\alpha_2 + A_{22}\alpha_1} $
		and $ 	A^2_{22}\alpha^2_1 + A^2_{21}\alpha^2_2 + (A_{21} - iA_{22})^2 = 0, $
		and $ \psi_{1}(s) $ satisfies 
		\begin{align*}
			\psi^{\prime}_{1}(s) + \psi_{1}(s+s_0) &= \left(\dfrac{A_{22}\alpha_1}{\alpha_1-\alpha_2}\right)\cos (L(z)+B) + \left(\dfrac{- A_{21}\alpha_2}{\alpha_1-\alpha_2}\right)\sin (L(z)+B)\\& -\left(\dfrac{-A_{22}}{\alpha_1-\alpha_2}\right)\sin(L(z)+L(c)+B) - \left(\dfrac{A_{21}}{\alpha_1-\alpha_2}\right)\cos(L(z)+L(c)+B).
		\end{align*}
	\end{enumerate}
\end{cor}
\begin{thm}\label{th-3.5}
Let $ a, b, C, \alpha, \beta, \gamma\in\mathbb{C} $ with $\alpha^2\neq ab $, $ \Delta\neq 0 $, $ c=(c_1,c_2)(\neq (0,0))\in\mathbb{C}^2 $, $ s_0=c_1+c_2 $ and $ f(z_1,z_2) $ be a finite order transcendental entire solutions for the general quadratic partial differential equation
\begin{align}\label{eq-3.5}
	aM_1(f)^2 + 2\alpha M_1(f)M_3(f)+ bM_3(f)^2 + 2\beta M_1(f) + 2\gamma M_3(f) + C  = 0.
\end{align}	
Then $ f(z_1,z_2) $ is one of the following forms
\begin{enumerate}
	\item [(i)] $ 	f(z_1,z_2)= A_0 + e^{(z_1 + Az_2 + B)}, $
     where $ K_9, K_{10}, B $ are constant with $ K^2_9 + K^2_{10}= 1 $, where
     \begin{align*}
     	\begin{cases}
     A_0 = \dfrac{K_9\xi^{\pm}_1}{\sqrt{\frac{DA^{\pm}}{-\Delta}}} - \dfrac{K_{10}\eta^{\pm}_1}{\sqrt{\frac{DB^{\mp}}{-\Delta}}} + \dfrac{\alpha\gamma- b\beta}{ab- \alpha^2},\; A=\dfrac{2k\pi i\pm\pi i-c_1}{c_2}, K_{10} =\dfrac{-R_2R_3\pm R_1\sqrt{R^2_1 +R^2_2 -R^2_3}}{R^2_1 + R^2_2},\vspace{2mm}\\
     R_1 =\dfrac{\xi^{\pm}_1-\eta^{\pm}_1}{\sqrt{\frac{DA^{\pm}}{-\Delta}}},\; R_2=\dfrac{\xi^{\pm}_1+ \eta^{\pm}_1}{\sqrt{\frac{DB^{\mp}}{-\Delta}}}\; \mbox{and}\; R_3=\dfrac{\beta(\alpha + b)- \gamma(\alpha + a)}{ab-\alpha^2}.
     	\end{cases}
     \end{align*}
     \item[(ii)] \begin{align*}
     	f(z_1,z_2)=(-e^{-c_1})^{\frac{z_2}{c_2}}\pi(z_2)e^{z_1} + R_{11}\sin [L(z)-L(c)+B] + R_{12}\cos [L(z)-L(c)+B],
     \end{align*}
     where $ L(z)=\alpha_1z_1 + \alpha_2z_2 $, $ \alpha_1,\alpha_2, B\in\mathbb{C} $ and
     \begin{align*}
  \begin{cases}
  e^{2iL(c)} = \dfrac{(iD_{11}-D_{12})\alpha_1 - (E_{11}-iE_{12})}{(iD_{11}+ D_{12})\alpha_1 + (iE_{12}+E_{11})},\vspace{2mm}\\
  (D^2_{11} + D^2_{12})\alpha^4_1 + 2(D_{11}E_{12} +D_{12}E_{11})\alpha^3_1 + (E^2_{11} + E^2_{12})\alpha^2_1 \\ \quad+ [(iD_{11} + D_{12})-(iE_{11}-E_{12})]=0,\; \mbox{where}\; \pi(z_2+c_2)=\pi(z_2).
   \end{cases}	
     \end{align*}   
\end{enumerate}
\end{thm}
The existence and form of trinomial quadratic PDEs is shown by the following result which is a consequence of Theorem \ref{th-3.5}.
\begin{cor}\label{cor-3.4}
Let $ a, b, C, \alpha\in\mathbb{C} $ with $\alpha^2\neq ab $, $ \Delta\neq 0 $, $ c=(c_1,c_2)(\neq (0,0))\in\mathbb{C}^2 $, $ s_0=c_1+c_2 $ and $ f(z_1,z_2) $ be a finite order transcendental entire solutions for the general quadratic partial differential equation
\begin{align*}
	a\left(f(z+c) + \frac{\partial f(z)}{\partial z_1}\right)^2 + 2\alpha \left(f(z+c) + \frac{\partial f(z)}{\partial z_1}\right)\left(f(z+c) + \dfrac{\partial^2 f(z)}{\partial z^2_1}\right)+ b\left(f(z+c) + \dfrac{\partial^2 f(z)}{\partial z^2_1}\right)^2 = 1.
\end{align*}	
Then $ f(z_1,z_2) $ is one of the following forms
\begin{enumerate}
	\item [(i)] $ f(z_1,z_2)= A_{31} + e^{(z_1 + Az_2 + B)}, $
	where $ K_{9*}, K_{10*}, B $ are constants with $ K^2_{9*} + K^2_{10*}= 1 $, where
	\[
	\begin{cases}
		A_{31} = A_{11}K_{9*} - A_{12}K_{10*},\; A=\dfrac{2k\pi i\pm\pi i-c_1}{c_2}\;\mbox{and}\; K_{10*} =\pm\dfrac{S_{11}}{\sqrt{S^2_{11} + S^2_{12}}},\vspace{2mm} \\
			S_{11} =\dfrac{\sqrt{2}(2\alpha - K_{14})}{K_{11}K_{12}}\; \mbox{and}\; S_{12} =\dfrac{\sqrt{2}(2\alpha + K_{14})}{K_{11}K_{12}}.
	\end{cases}
	\]
	\item[(ii)] 
	\begin{align*}
		f(z_1,z_2)=(-e^{-c_1})^{\frac{z_2}{c_2}}\pi(z_2)e^{z_1} + S_{15}\sin [L(z)-L(c)+B] + S_{16}\cos [L(z)-L(c)+B].
	\end{align*}
	where $ L(z)=\alpha_1z_1 + \alpha_2z_2 $, $ \alpha_1,\alpha_2, B\in\mathbb{C} $ and
	\begin{align*}
		e^{2iL(c)} = \dfrac{(iA_{11}-A_{12})\alpha_1 - (B_{11}-iB_{12})}{(iA_{11}+ A_{12})\alpha_1 + (iB_{12}+B_{11})},
	\end{align*}
	\begin{align*}
		(A^2_{11} + A^2_{12})\alpha^4_1 + &2(A_{11}B_{12} +A_{12}B_{11})\alpha^3_1 + (B^2_{11} + B^2_{12})\alpha^2_1 \\& + [(iA_{11} + A_{12})-(iB_{11}-B_{12})]=0,\; \mbox{where}\; \pi(z_2+c_2)=\pi(z_2)
	\end{align*}
	and $ S_{15}=\frac{(A_{11}\alpha_1 + B_{12}) - \alpha_1(A_{12}\alpha_1 + B_{11})}{\alpha^2_1 + 1}$, $S_{16}= \frac{(A_{12}\alpha_1 + B_{11}) + \alpha_1(A_{11}\alpha_1 + B_{12})}{\alpha^2_1 + 1}. $
\end{enumerate}
\end{cor}
We obtain the following corollary of Theorem \ref{th-3.5} which generalizes Theorem \ref{th-1.6}.
\begin{cor}\label{cor-3.8}
	Let $ a, b\in\mathbb{C} $ with $ \Delta\neq 0 $, $ c=(c_1,c_2)(\neq (0,0))\in\mathbb{C}^2 $, $ s_0=c_1+c_2 $ and $ f(z_1,z_2) $ be a finite order transcendental entire solutions for the general quadratic partial differential equation
	\begin{align*}
	a\left(f(z+c) + \frac{\partial f(z)}{\partial z_1}\right)^2 + b\left(f(z+c) + \dfrac{\partial^2 f(z)}{\partial z^2_1}\right)^2 = 1.
	\end{align*}	
	Then $ f(z_1,z_2) $ is one of the following forms
	\begin{enumerate}
		\item [(i)] $ f(z_1,z_2)= A_0 + e^{(z_1 + Az_2 + B)}, $	where $ K_{9**}, K_{10**}, B $ are constant with $ K^2_{9**} + K^2_{10**}= 1 $, $ A_0= -A_{21}K_{10**} $, $ A=\dfrac{2k\pi i\pm\pi i-c_1}{c_2} $, $ K_{10**} =\mp\frac{K_{21}}{\sqrt{K^2_{21}+K^2_{22}}}. $		
		\item[(ii)] 
		\begin{align*}
			f(z_1,z_2)=(-e^{-c_1})^{\frac{z_2}{c_2}}\pi(z_2)e^{z_1} + S_{17}\sin [L(z)-L(c)+B] + S_{18}\cos [L(z)-L(c)+B].
		\end{align*}
		where $ L(z)=\alpha_1z_1 + \alpha_2z_2 $, $ \alpha_1,\alpha_2, B\in\mathbb{C} $, $ e^{2iL(c)} = - 1 $ and $ A^2_{21}\alpha^4_1 + 2A_{21}A_{22}\alpha^3_1 + A^2_{22}\alpha^2_1 + (A_{21}-iA_{22})^2=0, $	where $ \pi(z_2+c_2)=\pi(z_2) $ and $ 	S_{17}=\frac{-\alpha_1(A_{21}\alpha_1 + A_{22})}{\alpha^2_1 + 1}, $ $S_{18}= \frac{(A_{21}\alpha_1 + A_{22})}{\alpha^2_1 + 1}. $
	\end{enumerate}
\end{cor}
\section{Proof of the main results}
First, we recall here some necessary lemmas which will play a key roles in proving the main results.
 	
\begin{lem}\cite{Ronkin_AMS-1974,Stoll-AMS-1974}\label{lem-4.1}
For any entire function $ F $ on $ \mathbb{C}^n $, $ F(0)\neq 0 $ and put $\rho(n_F)=\rho < \infty $. Then there exist a canonical function $ f_F $ and a function $ g_F \in\mathbb{C}^n $ such that $ F(z) =f_F (z)e^{g_F (z)} $. For the special case $ n = 1 $, $ f_F $ is the canonical product of Weierstrass.
 \end{lem}
\begin{rem}\label{rem-4.1}
Here, denote $ \rho(n_F) $ to be the order of the counting function of zeros of $ F $.
\end{rem}
\begin{lem}\cite{Polya-JLMS-1926}\label{lem-4.2}
If $  g $ and $ h $ are entire functions on the complex plane $ \mathbb{C} $ and $ g(h) $ is an entire function of finite order, then there are only two possible cases: either
\begin{enumerate}
 	\item [(i)] the internal function $ h $ is a polynomial and the external function $ g $ is of finite order; or
   \item[(ii)] the internal function $ h $ is not a polynomial but a function of finite order, and the external function $ g $ is of zero order.
\end{enumerate}
\end{lem}
\begin{lem}\cite{Hu-Li-Yang-2003}\label{lem-4.3}
Let $ f_j(\not \equiv 0) $, $ j=1,2,3 $, be meromorphic functions on $ \mathbb{C}^m $ such that $ f_1 $ is not constant and $ f_1+f_2+f_3=1 $ such that 
\begin{align*}
	\sum_{j=1}^{3}\left\{N_2(r,\frac{1}{f_j})+2\overline{N}(r,f_j)\right\}<\lambda T(r,f_1) + O(\log^{+} T(r,f_1)),
\end{align*}
for all $ r $ outside possibly a set with finite logarithmic measure, where $ \lambda<1 $ is a positive number. Then either $ f_2=1\;\mbox{or}\; f_3=1 $.
\end{lem}

Now we discuss the proof of the main results of the paper.
\begin{proof}[\bf Proof of Theorem \ref{th-3.1}]
Suppose that $ f $ is a transcendental entire solution of \eqref{eq-3.1} with finite order. Under the transformations discussed in Section 2, we that the equation \eqref{eq-3.1} can be written as
\begin{align}\label{eq-4.1}
	\left(\dfrac{DA^{\pm}}{-\Delta}\right)u^2 + \left(\dfrac{DB^{\mp}}{-\Delta}\right)v^2 = 1
\end{align}
 We divide the proof into two cases.

\noindent{\bf Case 1}: If $ \sqrt{\frac{DA^{\pm}}{-\Delta}}u $ is a constant, then it follows from \eqref{eq-4.1} that $ \sqrt{\frac{DB^{\mp}}{-\Delta}}v $ is also constant. Denoting
\begin{align*}
	\sqrt{\dfrac{DA^{\pm}}{-\Delta}}u = K_1\; \mbox{and}\; \sqrt{\dfrac{DB^{\mp}}{-\Delta}}v =K_2,
\end{align*}
in view of \eqref{eq-4.1}, it is easy to see that $ K^2_1 + K^2_2 = 1 $. This clearly forces to
\begin{align}\label{eq-4.2}
     L_1(f)= f(z) +  \dfrac{\partial f}{\partial z_1}=\dfrac{K_1\xi^{\pm}_1}{\sqrt{\frac{DA^{\pm}}{-\Delta}}} - \dfrac{K_2\eta^{\pm}_1}{\sqrt{\frac{DB^{\mp}}{-\Delta}}} + \dfrac{\alpha\gamma- b\beta}{ab- \alpha^2}:= A_1,\\ \label{eq-4.3}	
     L_2(f)=  f(z) +  \dfrac{\partial f}{\partial z_2}=\dfrac{K_1\eta^{\pm}_1}{\sqrt{\frac{DA^{\pm}}{-\Delta}}}+ \dfrac{K_2\xi^{\pm}_1}{\sqrt{\frac{DB^{\mp}}{-\Delta}}} + \dfrac{\alpha\beta- a\gamma}{ab- \alpha^2}:= B_1.
\end{align}
In view of \eqref{eq-4.2} and \eqref{eq-4.3}, we see that
\begin{align}\label{eq-4.4}
	\dfrac{\partial f}{\partial z_1} - \dfrac{\partial f}{\partial z_2} = A_1- B_1
\end{align}
The characteristic equations of \eqref{eq-4.4} are $ \dfrac{dz_1}{dt}= 1,\; \dfrac{d z_2}{dt} = -1,\; \dfrac{df}{dt}= A_1 -B_1. $
Using the initial conditions: $ z_1 = 0 $, $ z_2 = s $ and we thus get $ f =f(0,s):=\phi(s) $ with a parameter $ s $. Therefore, a simple computation shows that $ 	z_1 = t,\; z_2 = -t+s, $ and
\begin{align*}
	f(t,s)=\int_{0}^{t}(A_1 - B_1)dt + \phi(s)= (A_1 - B_1)t + \phi(s),
\end{align*}
where $ \phi(s) $ is a transcendental entire function with finite order in $ s $. Then we obtain $ z_1 =t $ and $ s= z_1 + z_2 $, the solution of the equation \eqref{eq-4.4} is of the form
\begin{align}\label{eq-4.5}
	f(z_1,z_2)= (A_1-B_1)z_1 + \phi(z_1 +z_2).
\end{align}
Since $ f $ is an entire function, we have $ \frac{\partial^2 f}{\partial z_1\partial z_2} = \frac{\partial^2 f}{\partial z_2\partial z_1} $, hence from \eqref{eq-4.2} and \eqref{eq-4.3}, we obtain $ \frac{\partial f}{\partial z_1} = \frac{\partial f}{\partial z_2} $. Therefore from \eqref{eq-4.4} it follows that $ A_1 = B_1 $. Furthermore, from \eqref{eq-4.2}, \eqref{eq-4.3} and using $ A_1 = B_1 $, $ K^2_1 + K^2_2 = 1 $, a simple computation shows that 
\begin{align*}
	K_2 =\dfrac{-R_2R_3\pm R_1\sqrt{R^2_1 +R^2_2 -R^2_3}}{R^2_1 + R^2_2}
\end{align*}
where 
\begin{align*}
R_1 =\dfrac{\xi^{\pm}_1-\eta^{\pm}_1}{\sqrt{\frac{DA^{\pm}}{-\Delta}}},\; R_2=\dfrac{\xi^{\pm}_1+ \eta^{\pm}_1}{\sqrt{\frac{DB^{\mp}}{-\Delta}}}\; \mbox{and}\; R_3=\dfrac{\beta(\alpha + b)- \gamma(\alpha + a)}{ab-\alpha^2}.
\end{align*}
Thus, the equation \eqref{eq-4.5} becomes
\begin{align}\label{eq-4.6}
	f(z_1,z_2)=\phi(z_1 +z_2).
\end{align}
Substituting \eqref{eq-4.6} into \eqref{eq-4.3}, we obtain $ 	\phi^{\prime}(z_1 + z_2) + \phi(z_1 + z_2) = A_1 $ and an easy computation shows that
\begin{align}\label{eq-4.7}
	\phi(z_1 + z_2)= A_1 + R_4 e^{-(z_1 +z_2)},\; R_4\in\mathbb{C}
\end{align}
Hence, from \eqref{eq-4.6} and \eqref{eq-4.7}, we see that the precise form of the solution is $ 	f(z_1,z_2)= A_1 + R_4 e^{-(z_1 +z_2)}. $

\noindent{\bf Case 2}:  If $ \sqrt{\frac{DA^{\pm}}{-\Delta}}u $ is not a constant, then in view of the fact that the entire solutions of the functional equation $ f^2 + g^2 = 1 $ are $ f= \cos h(z) $, $ g=\sin h(z) $, where $ h(z) $ is an entire function, it follows from \eqref{eq-4.1} and .... 
\begin{align}\label{eq-4.8}
	 L_1(f)=f(z) + \dfrac{\partial f}{\partial z_1}=D_{11}\cos h(z) - D_{12}\sin h(z) + T_1, \\\label{eq-4.9} L_2(f)= f(z) +  \dfrac{\partial f}{\partial z_2}=E_{11}\cos h(z) + E_{12}\sin h(z) + T_2.
\end{align}
 Subtracting \eqref{eq-4.9} from \eqref{eq-4.8}, we obtain 
\begin{align}\label{eq-4.10}
	\dfrac{\partial f}{\partial z_1} - \dfrac{\partial f}{\partial z_2}= (D_{11} - E_{11})\cos h(z) - (D_{12} + E_{12})\sin h(z) +(T_1 - T_2).
\end{align}
Using the fact that $ \frac{\partial^2 f}{\partial z_1\partial z_2} = \frac{\partial^2 f}{\partial z_2\partial z_1} $, from \eqref{eq-4.8} and \eqref{eq-4.9}, a routine computation shows that
\begin{align}\label{eq-4.11}
	\dfrac{\partial f}{\partial z_1} - \dfrac{\partial f}{\partial z_2}= \left(D_{11}\dfrac{\partial h}{\partial z_2} - E_{11}\dfrac{\partial h}{\partial z_1}\right)\sin h(z) + \left(D_{12}\dfrac{\partial h}{\partial z_2} + E_{12}\dfrac{\partial h}{\partial z_1}\right)\cos h(z)
\end{align}
From \eqref{eq-4.10} and \eqref{eq-4.11}, we obtain 
\begin{align}\label{eq-4.12}
	&\left(D_{11}\dfrac{\partial h}{\partial z_2} - E_{11}\dfrac{\partial h}{\partial z_1}+ (D_{12} + E_{12})\right)\sin h(z) \\&\nonumber + \left(D_{12}\dfrac{\partial h}{\partial z_2} + E_{12}\dfrac{\partial h}{\partial z_1}- (D_{11} - E_{11})\right)\cos h(z) = (T_1- T_2)
\end{align}
\noindent{\bf Sub-case 2.1}: Suppose that
\begin{align*}
	&D_{11}\dfrac{\partial h}{\partial z_2} - E_{11}\dfrac{\partial h}{\partial z_1}+ (D_{12} + E_{12})\neq 0, \\&D_{12}\dfrac{\partial h}{\partial z_2} + E_{12}\dfrac{\partial h}{\partial z_1}- (D_{11} - E_{11})\neq 0\; \mbox{and}\; (T_1- T_2)\neq 0.
\end{align*}
It follows from \eqref{eq-4.12} that
\begin{align*}
	\tan h(z) = -\dfrac{\left(D_{12}\dfrac{\partial h}{\partial z_2} + E_{12}\dfrac{\partial h}{\partial z_1}- (D_{11} - E_{11})\right)}{\left(D_{11}\dfrac{\partial h}{\partial z_2} - E_{11}\dfrac{\partial h}{\partial z_1}+ (D_{12} + E_{12})\right)} + \dfrac{(T_1- T_2)}{\cos h(z)\left(D_{11}\dfrac{\partial h}{\partial z_2} - E_{11}\dfrac{\partial h}{\partial z_1}+ (D_{12} + E_{12})\right)}.
\end{align*}
Since $ f $ is a finite order transcendental entire solution of the equation \eqref{eq-3.1}, by Lemmas \ref{lem-4.1}, \ref{lem-4.2}, \eqref{eq-4.8} and \eqref{eq-4.9}, we conclude that $ h(z) $ is a non-constant polynomial in $ \mathbb{C}^2 $. As $ h(z) $ is a non-constant polynomial, we assume that $ h(z)=\sum_{n=0}^{\infty}a_n z^n $. Using Nevanlinna theory we deduce that $ T(r, \tan h)=O\{T(r,h)+\log r + r^n\} $. But $ \lim_{r\rightarrow\infty}\frac{T(r, \tan h)}{T(r,h)+\log r + r^n}=\infty $ when $ h $ is a non-constant polynomial. Therefore $ h $ must be a constant, a contradiction.\vspace{1.2mm}

\noindent{\bf Sub-case 2.2}: Suppose that
\begin{align*}
	D_{11}\dfrac{\partial h}{\partial z_2} - E_{11}\dfrac{\partial h}{\partial z_1}+ (D_{12} + E_{12})= 0, D_{12}\dfrac{\partial h}{\partial z_2} + E_{12}\dfrac{\partial h}{\partial z_1}- (D_{11} - E_{11})\neq 0 , T_1- T_2\neq 0.
\end{align*}
Thus the equation \eqref{eq-4.12} can be written as 
\begin{align*}
	\tan h(z) =\dfrac{\left(D_{12}\dfrac{\partial h}{\partial z_2} + E_{12}\dfrac{\partial h}{\partial z_1}- (D_{11} - E_{11})\right)\sin h(z)}{T_1 - T_2}.
\end{align*}
Since $ f $ is a finite order transcendental entire solution of the equation \eqref{eq-3.1}, by Lemmas \ref{lem-4.1}, \ref{lem-4.2}, \eqref{eq-4.8} and \eqref{eq-4.9}, we conclude that $ h(z) $ is a non-constant polynomial in $ \mathbb{C}^2 $. Similarly as in the sub-case $ 2.1 $, we obtain that  $ \lim_{r\rightarrow\infty}\frac{T(r, \tan h)}{T(r,h)+\log r + r^n}=\infty $ when $ h $ is a non-constant polynomial. Therefore $ h $ must be a constant, a contradiction.\vspace{1.2mm}

\noindent{\bf Sub-case 2.3}: Assume 
\begin{align*}
D_{11}\dfrac{\partial h}{\partial z_2} - E_{11}\dfrac{\partial h}{\partial z_1}+ (D_{12} + E_{12})\neq 0, D_{12}\dfrac{\partial h}{\partial z_2} + E_{12}\dfrac{\partial h}{\partial z_1}- (D_{11} - E_{11})= 0, T_1- T_2\neq 0,
\end{align*}
By the similar argument used in Sub-case $ 1 $, we arrive a contradiction.

\noindent{\bf Sub-case 2.4}: Suppose that
\begin{align*}
D_{11}\dfrac{\partial h}{\partial z_2} - E_{11}\dfrac{\partial h}{\partial z_1}+ (D_{12} + E_{12})\neq 0, D_{12}\dfrac{\partial h}{\partial z_2} + E_{12}\dfrac{\partial h}{\partial z_1}- (D_{11} - E_{11})\neq 0, T_1- T_2= 0.
\end{align*}
By the similar argument used in Sub-case $ 1 $, we arrive a contradiction.

\noindent{\bf Sub-case 2.5}: Suppose that 
\begin{align*}
	D_{11}\dfrac{\partial h}{\partial z_2} - E_{11}\dfrac{\partial h}{\partial z_1}+ (D_{12} + E_{12})= 0, D_{12}\dfrac{\partial h}{\partial z_2} + E_{12}\dfrac{\partial h}{\partial z_1}- (D_{11} - E_{11})= 0, \;T_1- T_2\neq 0.
\end{align*}
In view of \eqref{eq-4.12}, we see that $ T_1- T_2= 0 $ which contradicts $ T_1- T_2\neq 0 $.

\noindent{\bf Sub-case 2.6}: Assume that
\begin{align*}
	D_{11}\dfrac{\partial h}{\partial z_2} - E_{11}\dfrac{\partial h}{\partial z_1}+ (D_{12} + E_{12})= 0, D_{12}\dfrac{\partial h}{\partial z_2} + E_{12}\dfrac{\partial h}{\partial z_1}- (D_{11} - E_{11})\neq 0,\;T_1- T_2= 0.
\end{align*}
Using \eqref{eq-4.12}, we obtain $ \cos h(z) =0 $ and this shows that $ h(z) $ is a constant, a contradiction. 

\noindent{\bf Sub-case 2.7}: Let
\begin{align*}
	D_{11}\dfrac{\partial h}{\partial z_2} - E_{11}\dfrac{\partial h}{\partial z_1}+ (D_{12} + E_{12})\neq 0, D_{12}\dfrac{\partial h}{\partial z_2} + E_{12}\dfrac{\partial h}{\partial z_1}- (D_{11} - E_{11})= 0, \;T_1- T_2= 0.
\end{align*}
In view of \eqref{eq-4.12} it is easy to see that $ \sin h(z) =0 $, this implies that $ h(z) $ is a constant, a contradiction. 

\noindent{\bf Sub-case 2.8}: Suppose that
\begin{align*}
D_{11}\dfrac{\partial h}{\partial z_2} - E_{11}\dfrac{\partial h}{\partial z_1}+ (D_{12} + E_{12})= 0, D_{12}\dfrac{\partial h}{\partial z_2} + E_{12}\dfrac{\partial h}{\partial z_1}- (D_{11} - E_{11})= 0, T_1- T_2= 0.
\end{align*}
Since $ T_1 = T_2 $, a simple computation shows that 
\begin{align}\label{eq-4.13}
	R_{11}\dfrac{\partial h}{\partial z_1} -R_{12}\dfrac{\partial h}{\partial z_2} = 0
\end{align}
where $ R_{11}= E_{11}(D_{11}-E_{11}) - E_{12}(D_{12}+E_{12}),\; R_{12}= D_{11}(D_{11}-E_{11}) + D_{12}(D_{12}+E_{12}). $ Thus, it follows from \eqref{eq-4.13} that
\begin{align}\label{eq-4.14}
	h(z)= z_2 +\dfrac{R_{12}}{R_{11}} z_1 + R_5,\;\; \mbox{where}\; R_5\in\mathbb{C}.
\end{align}
Substituting this $ h(z) $ into \eqref{eq-4.10}, we obtain
\begin{align}\label{eq-4.15}
	\dfrac{\partial f}{\partial z_1} - \dfrac{\partial f}{\partial z_2}&= (D_{11} - E_{11})\cos\left(z_2 +\dfrac{R_{12}}{R_{11}} z_1 + R_5\right) - (D_{12} + E_{12})\sin\left(z_2 +\dfrac{R_{12}}{R_{11}} z_1 + R_5\right).
\end{align}
The characteristic equations of \eqref{eq-4.15} are $ \dfrac{dz_1}{dt}= 1,\; \dfrac{d z_2}{dt} = -1 $, and hence
\begin{align*}
	\dfrac{df}{dt}=(D_{11} - E_{11})\cos\left(z_2 +\dfrac{R_{12}}{R_{11}} z_1 + R_5\right) - (D_{12} + E_{12})\sin\left(z_2 +\dfrac{R_{12}}{R_{11}} z_1 + R_5\right).
\end{align*}
Using the initial conditions: $ z_1=0 $, $ z_2=s $ and $ f =f(0,s):=\phi_1(s) $ with a parameter $ s $, wee see that $ z_1 = t,\; z_2 = -t+s, $ and
\begin{align*}
	f(t,s)&=\phi_1(s) \\&+ \int_{0}^{t} \left((D_{11} - E_{11})\cos\left(z_2 +\dfrac{R_{12}}{R_{11}} z_1 + R_5\right) -(D_{12} + E_{12})\sin\left(z_2 +\dfrac{R_{12}}{R_{11}} z_1 + R_5\right)\right)dt \\&= \phi_1(s) + (D_{11} - E_{11})\int_{0}^{t}\cos\left(\dfrac{R_{12}-R_{11}}{R_{11}}t + s + R_5\right)dt\\&\quad - (D_{12} + E_{12})\int_{0}^{t}\sin\left(\dfrac{R_{12}-R_{11}}{R_{11}}t + s + R_5\right)dt\\&= \phi_2(s) + \dfrac{(D_{11} - E_{11})R_{11}}{(R_{12}-R_{11})}\sin\left(\dfrac{R_{12}-R_{11}}{R_{11}}t + s + R_5\right) \\&\quad + \dfrac{(D_{12} + E_{12})R_{11}}{(R_{12}-R_{11})}\cos \left(\dfrac{R_{12}-R_{11}}{R_{11}}t + s + R_5\right)
\end{align*}
where $ \phi_2(s)=\phi_1(s) - \frac{(D_{11} - E_{11})R_{11}}{(R_{12}-R_{11})}\sin\left(s + R_5\right) - \frac{(D_{12} + E_{12})R_{11}}{(R_{12}-R_{11})}\cos \left(s + R_5\right) $ is a transcendental entire function with finite order in $ s $. Using $ t= z_1 $, $ s =z_1 +z_2 $, we easily obtain
\begin{align}\label{eq-4.16}
	f(z_1,z_2)&= \dfrac{(D_{11} - E_{11})R_{11}}{(R_{12}-R_{11})}\sin\left(\dfrac{R_{12}}{R_{11}}z_1 + z_2 + R_5\right) \\&\nonumber + \dfrac{(D_{12} + E_{12})R_{11}}{(R_{12}-R_{11})}\cos \left(\dfrac{R_{12}}{R_{11}}z_1 + z_2 + R_5\right) + \phi_2(z_1 + z_2)
\end{align}
Substituting \eqref{eq-4.16} into \eqref{eq-4.9} gives that
\begin{align*}
	\phi^{\prime}_2(z_1 + z_2) + \phi_2(z_1 + z_2) = R_{21}\sin\left(\dfrac{R_{12}}{R_{11}}z_1 + z_2 + R_5\right)+ R_{22}\cos\left(\dfrac{R_{12}}{R_{11}}z_1 + z_2 + R_5\right) + T_2,
\end{align*}
where
\begin{align*}
	R_{21}=\dfrac{E_{12}R_{12}+(E_{11} - D_{11} + D_{12})R_{11}}{(R_{12}-R_{11})}\;\mbox{and}\;R_{22}=\dfrac{E_{11}R_{12}-(E_{12} + D_{11} + D_{12})R_{11}}{(R_{12}-R_{11})}. 
\end{align*}
It follows easily that
\begin{align}\label{eq-4.17}
	\phi_2(z_1 + z_2)&=\dfrac{(R_{21}+ R_{22})}{2} \sin\left(\dfrac{R_{12}}{R_{11}}z_1 + z_2 + R_5\right) \\&\nonumber - \dfrac{(R_{21}- R_{22})}{2}\cos\left(\dfrac{R_{12}}{R_{11}}z_1 + z_2 + R_5\right) + R_6 e^{-(z_1 + z_2)} + T_2.
\end{align}
Therefore, from \eqref{eq-4.16} and \eqref{eq-4.17}, we obtain the precise form of the solution $ f $ as
\begin{align*}
	f(z_1,z_2)&=\left(\dfrac{E_{11}+E_{12}}{2}\right)\sin\left(\dfrac{R_{12}}{R_{11}}z_1 + z_2 + R_5\right) \\& + \left(\dfrac{E_{11}-E_{12}}{2}\right)\cos \left(\dfrac{R_{12}}{R_{11}}z_1 + z_2 + R_5\right) + R_6 e^{-(z_1 + z_2)} + T_2.
\end{align*}
This completes the proof.
\end{proof}

\begin{proof}[\bf Proof of Theorem \ref{th-3.2}]
Suppose that $ f $ is a transcendental entire solution of \eqref{eq-3.2} with finite order. We see that the equation \eqref{eq-3.2} can be written as
\begin{align}\label{eq-4.18}
	\left(\dfrac{DA^{\pm}}{-\Delta}\right)u^2 + \left(\dfrac{DB^{\mp}}{-\Delta}\right)v^2 = 1
\end{align}
We discuss the following two possible cases.

\noindent{\bf Case 1}: If $ \sqrt{\frac{DA^{\pm}}{-\Delta}}u $ is a constant, then it follows from \eqref{eq-4.18} that $ \sqrt{\frac{DB^{\mp}}{-\Delta}}v $ is also constant. Denote
\begin{align*}
	\sqrt{\dfrac{DA^{\pm}}{-\Delta}}u = K_3\; \mbox{and}\; \sqrt{\dfrac{DB^{\mp}}{-\Delta}}v =K_4.
\end{align*}
In view of \eqref{eq-4.18}, it is easy to see that $ K^2_3 + K^2_4 = 1 $, hence
\begin{align}\label{eq-4.19}
	L_1(f)= f(z) + \dfrac{\partial f}{\partial z_1}= \dfrac{K_3\xi^{\pm}_1}{\sqrt{\frac{DA^{\pm}}{-\Delta}}} - \dfrac{K_4\eta^{\pm}_1}{\sqrt{\frac{DB^{\mp}}{-\Delta}}} + \dfrac{\alpha\gamma- b\beta}{ab- \alpha^2}:= A_2,\\ \label{eq-4.20}	
	L_3(f) = f(z) + \dfrac{\partial^2 f}{\partial z^2_1}= \dfrac{K_1\eta^{\pm}_1}{\sqrt{\frac{DA^{\pm}}{-\Delta}}}+ \dfrac{K_2\xi^{\pm}_1}{\sqrt{\frac{DB^{\mp}}{-\Delta}}} + \dfrac{\alpha\beta- a\gamma}{ab- \alpha^2}:= B_2.
\end{align}
Subtracting \eqref{eq-4.20} from \eqref{eq-4.19}, we see that
\begin{align}\label{eq-4.21}
	\dfrac{\partial f}{\partial z_1} - \dfrac{\partial^2 f}{\partial z^2_1} =A_2 - B_2
\end{align}
Differentiating \eqref{eq-4.19} with respect to $ z_1 $ we obtain
\begin{align}\label{eq-4.22}
	\dfrac{\partial f}{\partial z_1} + \dfrac{\partial^2 f}{\partial z^2_1}=0
\end{align}
Therefore, it follows from \eqref{eq-4.20} and \eqref{eq-4.22} that
\begin{align}\label{eq-4.23}
	f(z) - \dfrac{\partial f}{\partial z_1} = B_2
\end{align}
Differentiating \eqref{eq-4.23} with respect to $ z_1 $ we obtain
\begin{align}\label{eq-4.24}
	\dfrac{\partial f}{\partial z_1} - \dfrac{\partial^2 f}{\partial z^2_1}= 0
\end{align}
In view of \eqref{eq-4.21} and \eqref{eq-4.24}, it is clear that $ A_2 = B_2 $. Using $ K^3_1 + K^4_2 = 1 $, $ A_2 = B_2 $ and from \eqref{eq-4.2} and \eqref{eq-4.3}, a simple computation shows that 
\begin{align*}
	K_3 =\dfrac{-R_2R_3\pm R_1\sqrt{R^2_1 +R^2_2 -R^2_3}}{R^2_1 + R^2_2}
\end{align*}
where 
\begin{align*}
	R_1 =\dfrac{\xi^{\pm}_1-\eta^{\pm}_1}{\sqrt{\frac{DA^{\pm}}{-\Delta}}},\; R_2=\dfrac{\xi^{\pm}_1+ \eta^{\pm}_1}{\sqrt{\frac{DB^{\mp}}{-\Delta}}}\; \mbox{and}\; R_3=\dfrac{\beta(\alpha + b)- \gamma(\alpha + a)}{ab-\alpha^2}.
\end{align*}
Therefore, from \eqref{eq-4.24} it follows that
\begin{align}\label{eq-4.25}
	f(z_1,z_2)= e^{z_1 + \phi_1(z_2)} + \phi_2(z_2),
\end{align}
where $ \phi_1(z_2) $, $ \phi_2(z_2) $ are two functions in $ z_2 $. Substituting \eqref{eq-4.25} into \eqref{eq-4.19} we obtain $ 	2 e^{z_1 + \phi_1(z_2)} + \phi_2(z_2) = A_2. $ Thus, it follows that $ \phi_2(z_2) = A_2 $ and $ 2 e^{z_1 + \phi_1(z_2)}=0 $, which is impossible.

\noindent{\bf Case 2}:  If $ \sqrt{\frac{DA^{\pm}}{-\Delta}}u $ is not a constant. Similar to the argument as in the proof of Theorem \ref{th-3.1}, there exists an entire function $ h(z) $ such that  
\begin{align}\label{eq-4.26}
	L_1(f)=f(z) + \dfrac{\partial f}{\partial z_1}=D_{11}\cos h(z) - D_{12}\sin h(z) + T_1, \\\label{eq-4.27} L_3(f)= f(z) + \dfrac{\partial^2 f}{\partial z^2_1}=E_{11}\cos h(z) + E_{12}\sin h(z) + T_2.
\end{align}
By the similar reasoning applied previously, we obtain 
\begin{align}\label{eq-4.28}
    \dfrac{\partial f}{\partial z_1} - \dfrac{\partial^2 f}{\partial z^2_1}= (D_{11} - E_{11})\cos h(z) - (D_{12} + E_{12})\sin h(z) +(T_1 - T_2).
\end{align}
Differentiating \eqref{eq-4.26} with respect to $ z_1 $ it follows that
\begin{align}\label{eq-4.29}
	\dfrac{\partial f}{\partial z_1} + \dfrac{\partial^2 f}{\partial z^2_1}= -D_{11}\sin h(z) \dfrac{\partial h}{\partial z_1} - D_{12}\cos h(z) \dfrac{\partial h}{\partial z_1}
\end{align}
From \eqref{eq-4.27} and \eqref{eq-4.29} it is easy to see that
\begin{align}\label{eq-4.30}
	f(z) - \dfrac{\partial f}{\partial z_1} =\left(E_{11} + D_{12}\dfrac{\partial h}{\partial z_1}\right)\cos h(z) + \left(E_{12} + D_{11}\dfrac{\partial h}{\partial z_1}\right)\sin h(z) + T_2.
\end{align}
Differentiating \eqref{eq-4.30} with respect to $ z_1 $ and using \eqref{eq-4.28} we obtain
\begin{align*}
	&\left((D_{12} + E_{12}) - E_{11}\dfrac{\partial h}{\partial z_1} - D_{12}\left(\dfrac{\partial h}{\partial z_1}\right)^2 + D_{11}\dfrac{\partial^2 h}{\partial z^2_1}\right)\sin h(z) \\&= \left((D_{11} - E_{11}) - E_{12}\dfrac{\partial h}{\partial z_1} - D_{11}\left(\dfrac{\partial h}{\partial z_1}\right)^2 - D_{12}\dfrac{\partial^2 h}{\partial z^2_1}\right)\cos h(z) + (T_1 - T_2)
\end{align*}
Since $ f $ is a finite order transcendental entire solution of the equation \eqref{eq-3.2}, by the Lemmas \ref{lem-4.1}, \ref{lem-4.2}, \eqref{eq-4.26} and \eqref{eq-4.27}, we conclude that $ h(z) $ is a non-constant polynomial in $ \mathbb{C}^2 $. Similar arguments as in the proof of Theorem \ref{th-3.1} gives that
\begin{align}\label{eq-4.31}
	&(D_{12} + E_{12}) - E_{11}\dfrac{\partial h}{\partial z_1} - D_{12}\left(\dfrac{\partial h}{\partial z_1}\right)^2 + D_{11}\dfrac{\partial^2 h}{\partial z^2_1} = 0,\\&\label{eq-4.32} (D_{11} - E_{11}) - E_{12}\dfrac{\partial h}{\partial z_1} - D_{11}\left(\dfrac{\partial h}{\partial z_1}\right)^2 - D_{12}\dfrac{\partial^2 h}{\partial z^2_1} = 0, \\&\label{eq-4.33}\mbox{and}\; (T_1 - T_2)= 0.
\end{align}
In view of \eqref{eq-4.31} and \eqref{eq-4.32}, we see that
\begin{align}\label{eq-4.34}
	\dfrac{\partial^2 h}{\partial z^2_1} - S_1\dfrac{\partial h}{\partial z_1} = S_2,
\end{align}
where $ S_1 =\frac{D_{11}E_{11}- D_{12}E_{12}}{D^2_{11}+D^2_{12}} \;\;\mbox{and}\;\; S_2 =-\frac{D_{11}E_{12}+ D_{12}E_{11}}{D^2_{11}+ D^2_{12}} $ which implies that $ h(z_1,z_2)=e^{S_1z_1+\phi_1(z_2)} + e^{\phi_2(z_2)} - \frac{S_2}{S_1}z_1  $, is a transcendental polynomial, which is a contradiction as $ h(z) $ is a non-constant polynomial. This completes the proof of Theorem \ref{th-3.2}.
\end{proof}	

\begin{proof}[\bf Proof of Theorem \ref{th-3.3}]
Suppose that $ f $ is a transcendental entire solution of \eqref{eq-3.3} with finite order. The equation \eqref{eq-3.3} can be written as
\begin{align}\label{eq-4.35}
	\left(\dfrac{DA^{\pm}}{-\Delta}\right)u^2 + \left(\dfrac{DB^{\mp}}{-\Delta}\right)v^2 = 1
\end{align}
Now we will discuss two cases below.

\noindent{\bf Case 1}: If $ \sqrt{\frac{DA^{\pm}}{-\Delta}}u $ is a constant, then it follows from \eqref{eq-4.35} that $ \sqrt{\frac{DB^{\mp}}{-\Delta}}v $ is also constant. Denote
\begin{align*}
	\sqrt{\dfrac{DA^{\pm}}{-\Delta}}u = K_5\;\mbox{and}\; \sqrt{\dfrac{DB^{\mp}}{-\Delta}}v =K_6
\end{align*}
In view of \eqref{eq-4.35}, it is easy to see that $ K^2_5 + K^2_6 = 1 $. Thus, we have
\begin{align}\label{eq-4.36}
	L_1(f)= f(z) + \dfrac{\partial f}{\partial z_1}=\dfrac{K_5\xi^{\pm}_1}{\sqrt{\frac{DA^{\pm}}{-\Delta}}} - \dfrac{K_6\eta^{\pm}_1}{\sqrt{\frac{DB^{\mp}}{-\Delta}}} + \dfrac{\alpha\gamma- b\beta}{ab- \alpha^2}:= A_3,\\ \label{eq-4.37}	
	L_4(f)= f(z) + \dfrac{\partial^2 f}{\partial z_1 \partial z_2} =\dfrac{K_5\eta^{\pm}_1}{\sqrt{\frac{DA^{\pm}}{-\Delta}}}+ \dfrac{K_6\xi^{\pm}_1}{\sqrt{\frac{DB^{\mp}}{-\Delta}}} + \dfrac{\alpha\beta- a\gamma}{ab- \alpha^2}:= B_3.
\end{align}
An easy computation using \eqref{eq-4.36} and \eqref{eq-4.37} shows that
\begin{align}\label{eq-4.38}
	\dfrac{\partial f}{\partial z_1} - \dfrac{\partial^2 f}{\partial z_1\partial z_2} =A_3 - B_3
\end{align}
Differentiating \eqref{eq-4.36} with respect to $ z_2 $ we obtain
\begin{align}\label{eq-4.39}
	\dfrac{\partial f}{\partial z_2} + \dfrac{\partial^2 f}{\partial z_2\partial z_1}=0
\end{align}
From \eqref{eq-4.37} and \eqref{eq-4.39}, it follows that
\begin{align}\label{eq-4.40}
	f(z) - \dfrac{\partial f}{\partial z_2} = B_3
\end{align}
Differentiating \eqref{eq-4.40} with respect to $ z_1 $, we obtain
\begin{align}\label{eq-4.41}
	\dfrac{\partial f}{\partial z_1} - \dfrac{\partial^2 f}{\partial z_1\partial z_2}= 0
\end{align}
In view of \eqref{eq-4.39} and \eqref{eq-4.41} it is easy to see that $ A_3 = B_3 $. Using $ K^2_5 + K^2_6 = 1 $, $ A_3 = B_3 $ and from \eqref{eq-4.36}, \eqref{eq-4.37} a simple computation shows that 
\begin{align*}
	K_6 =\dfrac{-R_2R_3\pm R_1\sqrt{R^2_1 +R^2_2 -R^2_3}}{R^2_1 + R^2_2}
\end{align*}
where 
\begin{align*}
	R_1 =\dfrac{\xi^{\pm}_1-\eta^{\pm}_1}{\sqrt{\frac{DA^{\pm}}{-\Delta}}},\; R_2=\dfrac{\xi^{\pm}_1+ \eta^{\pm}_1}{\sqrt{\frac{DB^{\mp}}{-\Delta}}}\; \mbox{and}\; R_3=\dfrac{\beta(\alpha + b)- \gamma(\alpha + a)}{ab-\alpha^2}.
\end{align*}
Combining \eqref{eq-4.39} and \eqref{eq-4.41} yields that
\begin{align}\label{eq-4.42}
	\dfrac{\partial f}{\partial z_1} + \dfrac{\partial f}{\partial z_2} =0,
\end{align}
which shows that $ f(z_1,z_2)=\phi(z_2 - z_1) $. Substituting $ f(z_1,z_2) $ into \eqref{eq-4.36}, we obtain $ \phi^{\prime}(z_2 - z_1) - \phi(z_2 - z_1) = -A_3. $ Of course, $ \phi(z_2 - z_1) = -A_3 + B_{11}e^{(z_1 - z_2)} $, $ B_{11}\in\mathbb{C} $. Hence, $ f $ takes the form $ f(z_1,z_2)= -A_3 + B_{11}e^{(z_1 - z_2)} $.

\noindent{\bf Case 2}:  If $ \sqrt{\frac{DA^{\pm}}{-\Delta}}u $ is not a constant, similar argument being used in the proof of Theorem \ref{th-3.1} that, there exists an entire function $ h(z) $ such that  
\begin{align}\label{eq-4.43}
	L_1(f)= f(z) + \dfrac{\partial f}{\partial z_1}=D_{11}\cos h(z) - D_{12}\sin h(z) + T_1, \\\label{eq-4.44} L_4(f)= f(z) + \dfrac{\partial^2 f}{\partial z_1 \partial z_2}=E_{11}\cos h(z) + E_{12}\sin h(z) + T_2.
\end{align}
In view of \eqref{eq-4.43} and \eqref{eq-4.44}, a simple computation shows that \begin{align}\label{eq-4.45}
	\dfrac{\partial f}{\partial z_1} - \dfrac{\partial^2 f}{\partial z_1 \partial z_2}= (D_{11} - E_{11})\cos h(z) - (D_{12} + E_{12})\sin h(z) +(T_1 - T_2).
\end{align}
Differentiating \eqref{eq-4.43} with respect to $ z_2 $, we obtain
\begin{align}\label{eq-4.46}
	\dfrac{\partial f}{\partial z_2} + \dfrac{\partial^2 f}{\partial z_2\partial z_1}= - D_{11}\sin h(z)\dfrac{\partial h}{\partial z_2} - D_{12}\cos h(z)\dfrac{\partial h}{\partial z_2}
\end{align}
From \eqref{eq-4.44} and \eqref{eq-4.46}, it follows that
\begin{align}\label{eq-4.47}
	f(z) - \dfrac{\partial f}{\partial z_2} =\left(D_{11}\dfrac{\partial h}{\partial z_2} + E_{12}\right)\sin h(z) + \left(D_{12}\dfrac{\partial h}{\partial z_2} + E_{11}\right)\cos h(z) + T_2
\end{align}
Differentiating \eqref{eq-4.47} with respect to $ z_1 $, and using \eqref{eq-4.45}, the following can be obtained
\begin{align}\label{eq-4.48}
	&\left(D_{11}\dfrac{\partial^2 h}{\partial z_1\partial z_2} - D_{12}\dfrac{\partial h}{\partial z_1}\dfrac{\partial h}{\partial z_2}- E_{11}\dfrac{\partial h}{\partial z_1} + (D_{12} + E_{12})\right)\sin h(z) \\&\nonumber = \left((D_{11}-E_{11}) - D_{12}\dfrac{\partial^2 h}{\partial z_1\partial z_2}- D_{11}\dfrac{\partial h}{\partial z_1}\dfrac{\partial h}{\partial z_2} - E_{12}\dfrac{\partial h}{\partial z_1}\right)\cos h(z) + (T_1 - T_2).
\end{align}
However, since $ f $ is a finite order transcendental entire solution of the equation \eqref{eq-3.3}, by Lemmas \ref{lem-4.1}, \ref{lem-4.2}, \eqref{eq-4.43} and \eqref{eq-4.44}, we conclude that $ h(z) $ is a non-constant polynomial in $ \mathbb{C}^2 $. The similar argument being used in the proof of Theorem \ref{th-3.1} will help us to obtain the following
\begin{align}\label{eq-4.49}
	&D_{11}\dfrac{\partial^2 h}{\partial z_1\partial z_2} - D_{12}\dfrac{\partial h}{\partial z_1}\dfrac{\partial h}{\partial z_2}- E_{11}\dfrac{\partial h}{\partial z_1} + (D_{12} + E_{12}) = 0,\\&\label{eq-4.50} (D_{11}-E_{11}) - D_{12}\dfrac{\partial^2 h}{\partial z_1\partial z_2}- D_{11}\dfrac{\partial h}{\partial z_1}\dfrac{\partial h}{\partial z_2} - E_{12}\dfrac{\partial h}{\partial z_1} = 0 \\&\label{eq-4.51}\mbox{and}\; (T_1 - T_2)= 0.
\end{align}
In view of \eqref{eq-4.49} and \eqref{eq-4.50}, a simple computation shows that $ 	\frac{\partial h}{\partial z_1}\frac{\partial h}{\partial z_2} + S_3 \frac{\partial h}{\partial z_1} = S_4, $ where 
\begin{align*}
	S_3= \dfrac{D_{11}E_{12} + D_{12}E_{11}}{D^2_{11} + D^2_{12}}\;\;\mbox{and}\;\; S_4=\dfrac{D_{11}(D_{11}-E_{11}) + D_{12}(D_{12} + E_{12})}{D^2_{11} + D^2_{12}}
\end{align*}
which implies that $ h(z_1,z_2)= \phi(z_1)e^{-S_3 z_2} + \phi(z_2) +\frac{S_4}{S_3}z_1 $, is a transcendental polynomial, which is a contradiction as $ h(z) $ is a non-constant polynomial. This completes the proof of Theorem \ref{th-3.3}.
\end{proof}

\begin{proof}[\bf Proof of Theorem \ref{th-3.4}]
Suppose that $ f $ is a transcendental entire solution of \eqref{eq-3.4} with finite order. Under the transformation mentioned in Section 2,  \eqref{eq-3.4} can be written as
\begin{align}\label{eq-4.52}
	\left(\dfrac{DA^{\pm}}{-\Delta}\right)u^2 + \left(\dfrac{DB^{\mp}}{-\Delta}\right)v^2 = 1
\end{align}

\noindent{\bf Case 1}: If $ \sqrt{\frac{DA^{\pm}}{-\Delta}}u $ is a constant, then it follows from \eqref{eq-4.52} that $ \sqrt{\frac{DB^{\mp}}{-\Delta}}v $ is also constant. Denote
\begin{align*}
	\sqrt{\dfrac{DA^{\pm}}{-\Delta}}u = K_7\; \mbox{and}\; \sqrt{\dfrac{DB^{\mp}}{-\Delta}}v =K_8.
\end{align*}
From \eqref{eq-4.52}, it is clear that $ K^2_7 + K^2_8 = 1 $. As a matter of fact, we see that 
\begin{align}\label{eq-4.53}
	M_1(f)= f(z+c) + \dfrac{\partial f}{\partial z_1}
	= \dfrac{K_7\xi^{\pm}_1}{\sqrt{\frac{DA^{\pm}}{-\Delta}}} - \dfrac{K_8\eta^{\pm}_1}{\sqrt{\frac{DB^{\mp}}{-\Delta}}} + \dfrac{\alpha\gamma- b\beta}{ab- \alpha^2}:= A_4,\\ \label{eq-4.54}	
	M_2(f) = f(z+c) + \dfrac{\partial f}{\partial z_2}= \dfrac{K_7\eta^{\pm}_1}{\sqrt{\frac{DA^{\pm}}{-\Delta}}}+ \dfrac{K_8\xi^{\pm}_1}{\sqrt{\frac{DB^{\mp}}{-\Delta}}} + \dfrac{\alpha\beta- a\gamma}{ab- \alpha^2}:= B_4.
\end{align}
From \eqref{eq-4.53} and \eqref{eq-4.54}, we obtain
\begin{align}\label{eq-4.55}
	\dfrac{\partial f}{\partial z_1} - \dfrac{\partial f}{\partial z_2} =A_4 - B_4.
\end{align}
The characteristic equations of \eqref{eq-4.55} are $ 	\frac{dz_1}{dt}= 1,\; \frac{d z_2}{dt} = -1,\; \frac{df}{dt}= A_4 -B_4. $ Using the initial conditions: $ z_1 = 0 $, $ z_2 = s $ and $ f =f(0,s):=\psi(s) $ with a parameter $ s $. Therefore, a simple computation shows that $ 	z_1 = t,\; z_2 = -t+s, $ and $ 	f(t,s) =\int_{0}^{t}(A_4 - B_4)dt + \psi(s)= (A_4 - B_4)t + \psi(s), $ where $ \psi(s) $ is a transcendental entire function with finite order in $ s $. Then we obtain $ z_1 =t $ and $ s= z_1 + z_2 $, the solution of the equation \eqref{eq-4.55} is of the form
\begin{align}\label{eq-4.56}
	f(z_1,z_2)= (A_4 -B_4)z_1 + \psi(z_1 +z_2).
\end{align}
In view of $ \frac{\partial^2 f}{\partial z_1\partial z_2} = \frac{\partial^2 f}{\partial z_2\partial z_1} $, using \eqref{eq-4.2} and \eqref{eq-4.3}, we obtain that $ \frac{\partial f(z+c)}{\partial z_1} = \frac{\partial f(z+c)}{\partial z_2} $. Thus, it follows from \eqref{eq-4.55} that $ A_4 = B_4 $. Also, from \eqref{eq-4.53}, \eqref{eq-4.54} and using $ A_0 = B_0 $, $ K^2_7 + K^2_8 = 1 $, a simple computation shows that $ 	K_8 =\frac{-R_2R_3\pm R_1\sqrt{R^2_1 +R^2_2 -R^2_3}}{R^2_1 + R^2_2} $, where 
\begin{align*}
	R_1 =\dfrac{\xi^{\pm}_1-\eta^{\pm}_1}{\sqrt{\frac{DA^{\pm}}{-\Delta}}},\; R_2=\dfrac{\xi^{\pm}_1+ \eta^{\pm}_1}{\sqrt{\frac{DB^{\mp}}{-\Delta}}}\; \mbox{and}\; R_3=\dfrac{\beta(\alpha + b)- \gamma(\alpha + a)}{ab-\alpha^2}.
\end{align*}
From \eqref{eq-4.56} and using $ A_4 = B_4 $ we obtain $ f(z_1,z_2)= \psi(z_1 +z_2) $, where $ \psi(s) $ is a transcendental entire function with finite order in $ s $. Also $ \psi(s) $ satisfying the relation $ \psi^{\prime}(s) + \psi(s+ s_0) = A_4 $, where $ s_0= c_1 +c_2 $.

\noindent{\bf Case 2}:  If $ \sqrt{\frac{DA^{\pm}}{-\Delta}}u $ is not a constant. By the similar argument being used previously, we obtain   
\begin{align}\label{eq-4.57}
M_1(f)= f(z+c) + \dfrac{\partial f}{\partial z_1}=D_{11}\cos h(z) - D_{12}\sin h(z) + T_1, \\\label{eq-4.58} M_2(f) = f(z+c) + \dfrac{\partial f}{\partial z_2}=E_{11}\cos h(z) + E_{12}\sin h(z) + T_2.
\end{align}
Subtracting \eqref{eq-4.58} from \eqref{eq-4.57}, we see that 
\begin{align}\label{eq-4.59}
\dfrac{\partial f}{\partial z_1} - \dfrac{\partial f}{\partial z_2}= (D_{11} - E_{11})\cos h(z) - (D_{12} + E_{12})\sin h(z) +(T_1 - T_2).
\end{align}
Since $ \frac{\partial^2 f}{\partial z_1\partial z_2} = \frac{\partial^2 f}{\partial z_2\partial z_1} $, using \eqref{eq-4.57} and \eqref{eq-4.58}, a simple computation shows that
\begin{align}\label{eq-4.60}
\dfrac{\partial f(z+c)}{\partial z_1} - \dfrac{\partial f(z+c)}{\partial z_2}= \left(D_{11}\dfrac{\partial h}{\partial z_2}- E_{11}\dfrac{\partial h}{\partial z_1}\right)\sin h(z) + \left(D_{12}\dfrac{\partial h}{\partial z_2}+ E_{12}\dfrac{\partial h}{\partial z_1}\right)\cos h(z).
\end{align}
Therefore, it follows from \eqref{eq-4.59} and \eqref{eq-4.60} that
\begin{align}\label{eq-4.61}
	&\left(\dfrac{\{(D_{12}-iD_{11})+(E_{12}+iE_{11})\}e^{ih(z+c)} + 2i(T_2 - T_1)}{(iD_{11} + D_{12})-(iE_{11}-E_{12})}\right)e^{ih(z+c)}\\&\nonumber + \left(\dfrac{(D_{11}+iD_{12})\dfrac{\partial h}{\partial z_2} -(E_{11}-iE_{12})\dfrac{\partial h}{\partial z_1}}{(iD_{11} + D_{12})-(iE_{11}-E_{12})}\right)e^{i[h(z+c)+h(z)]}\\&\nonumber + \left(\dfrac{(E_{11}+iE_{12})\dfrac{\partial h}{\partial z_1}-(D_{11}-iD_{12})\dfrac{\partial h}{\partial z_2}}{(iD_{11} + D_{12})-(iE_{11}-E_{12})}\right)e^{i[h(z+c)-h(z)]} = 1.
\end{align}
If $ \left(\frac{(D_{11}+iD_{12})\frac{\partial h}{\partial z_2} -(E_{11}-iE_{12})\frac{\partial h}{\partial z_1}}{(iD_{11} + D_{12})-(iE_{11}-E_{12})}\right) = 0 $, then we have $ \left(\frac{(E_{11}+iE_{12})\frac{\partial h}{\partial z_1}-(D_{11}-iD_{12})\frac{\partial h}{\partial z_2}}{(iD_{11} + D_{12})-(iE_{11}-E_{12})}\right)\neq 0 $. Otherwise, $ \left(\frac{\{(D_{12}-iD_{11})+(E_{12}+iE_{11})\}e^{ih(z+c)} + 2i(T_2 - T_1)}{(iD_{11} + D_{12})-(iE_{11}-E_{12})}\right)e^{ih(z+c)} = 1 $, thus, it follows that $ e^{ih(z+c)}=\frac{i(T_1 - T_2)\pm i\sqrt{(T_2-T_1)^2 + [(D_{12}-iD_{11})+(E_{12}+iE_{11})][(iE_{11}-E_{12})-(iD_{11}+D_{12})]}}{(D_{12}-iD_{11})+(E_{12}+iE_{11})} $, which implies that $ h(z) $ must be a constant, and hence, $ \sqrt{\frac{DA^{\pm}}{-\Delta}}u $ is a constant, which is a contradiction. Therefore, the equation \eqref{eq-4.61} is reduced to the following form
\begin{align*}
&\left(\dfrac{\{(D_{12}-iD_{11})+(E_{12}+iE_{11})\}e^{ih(z+c)} + 2i(T_2 - T_1)}{(iD_{11} + D_{12})-(iE_{11}-E_{12})}\right)e^{ih(z+c)}\\&\nonumber + \left(\dfrac{(E_{11}+iE_{12})\dfrac{\partial h}{\partial z_1}-(D_{11}-iD_{12})\dfrac{\partial h}{\partial z_2}}{(iD_{11} + D_{12})-(iE_{11}-E_{12})}\right)e^{i[h(z+c)-h(z)]}= 1.
\end{align*}
By the similar argument being used in \cite{Xu-Tu-Wang-RM-2021}, using Second Fundamental Theorem of Nevanlinna for several complex variables, we obtain $ T(r,e^{ih(z+c)})=o\left(T\left(r,e^{h}\right)\right) $, which is a contradiction. Therefore, $ \left(\frac{(D_{11}+iD_{12})\frac{\partial h}{\partial z_2} -(E_{11}-iE_{12})\frac{\partial h}{\partial z_1}}{(iD_{11} + D_{12})-(iE_{11}-E_{12})}\right)\neq 0 $. Similarly, it can be shown that  $ \left(\frac{(E_{11}+iE_{12})\frac{\partial h}{\partial z_1}-(D_{11}-iD_{12})\frac{\partial h}{\partial z_2}}{(iD_{11} + D_{12})-(iE_{11}-E_{12})}\right)\neq 0 $. Thus, in view of Lemma \ref{lem-4.3} and \eqref{eq-4.61}, we easily obtain
\begin{align}\label{eq-4.63}
	\left(\dfrac{(E_{11}+iE_{12})\dfrac{\partial h}{\partial z_1}-(D_{11}-iD_{12})\dfrac{\partial h}{\partial z_2}}{(iD_{11} + D_{12})-(iE_{11}-E_{12})}\right)e^{i[h(z+c)-h(z)]}\equiv 1
\end{align}
Also, from \eqref{eq-4.61}, it follows that
\begin{align}\label{eq-4.64}
	\left(\dfrac{(D_{11}+iD_{12})\dfrac{\partial h}{\partial z_2} -(E_{11}-iE_{12})\dfrac{\partial h}{\partial z_1}}{(iD_{11} + D_{12})-(iE_{11}-E_{12})}\right)e^{i[h(z)-h(z+c)]}\equiv 1.
\end{align}
Since $ h(z) $ is a polynomial, from \eqref{eq-4.63}(or \eqref{eq-4.64}), it follows that $ h(z+c)-h(z)=\xi $, where $ \xi $ is a constant in $ \mathbb{C} $. Hence, we see that $ h(z)=L(z) + H(z) + B $, where $ L(z)=\alpha_1z_1 + \alpha_2z_2 $, $ H(z):= H(s_1) $, $ H(s_1) $ is a polynomial in $ s_1=c_2z_1 - c_1z_2 $. We show that $ H(z)\equiv 0 $. Using \eqref{eq-4.63} and \eqref{eq-4.64}, a simple computation shows that 
\begin{align}\label{eq-4.65}
	\dfrac{(E_{11}+iE_{12})\alpha_1 - (D_{11}-iD_{12})\alpha_2 + [(D_{11}-iD_{12})c_1 + (E_{11}+iE_{12})c_2]H^{\prime}}{(iD_{11} + D_{12})-(iE_{11}-E_{12})}= e^{-iL(c)}
\end{align}
and
\begin{align}\label{eq-4.66}
	\dfrac{(D_{11}+iD_{12})\alpha_2 - (E_{11}-iE_{12})\alpha_1 - [(D_{11}+iD_{12})c_1 + (E_{11}-iE_{12})c_2]H^{\prime}}{(iD_{11} + D_{12})-(iE_{11}-E_{12})}= e^{iL(c)}.
\end{align}
In view of \eqref{eq-4.65} and \eqref{eq-4.66}, now it follows that $ [(D_{11}-iD_{12})c_1 + (E_{11}+iE_{12})c_2]H^{\prime} $ and $ [(D_{11}+iD_{12})c_1 + (E_{11}-iE_{12})c_2]H^{\prime} $ are constants. Since $ (c_1,c_2)\neq (0,0) $, it follows that $ H^{\prime} $ is a constant, that is $ \deg_{s_1}H\leq 1 $. Hence, it is easy to see that $ L(z) + H(z) + B $ is still linear form of $ \alpha_1z_1 + \alpha_2z_2 + B $, which means that $ H(z)\equiv 0 $. Therefore it is easy to see that $ h(z)=L(z) + B =\alpha_1z_1 + \alpha_2z_2 + B $. Now, in view of \eqref{eq-4.65} and \eqref{eq-4.66}, we see that
\[
\begin{cases}
		\dfrac{(E_{11}+iE_{12})\alpha_1 - (D_{11}-iD_{12})\alpha_2}{(iD_{11} + D_{12})-(iE_{11}-E_{12})} e^{iL(c)} = 1;\vspace{2mm}\\\dfrac{(D_{11}+iD_{12})\alpha_2 - (E_{11}-iE_{12})\alpha_1}{(iD_{11} + D_{12})-(iE_{11}-E_{12})}e^{-iL(c)}= 1,
		
\end{cases}
\]
which implies that
\[
\begin{cases}
e^{2iL(c)}=\dfrac{(D_{11}+iD_{12})\alpha_2 - (E_{11}-iE_{12})\alpha_1}{(E_{11}+iE_{12})\alpha_1 - (D_{11}-iD_{12})\alpha_2},\vspace{2mm}\\
(E^2_{11} + E^2_{12})\alpha^2_1 + (D^2_{11} + D^2_{12})\alpha^2_2 - 2(D_{11}E_{11} - D_{12}E_{12})\alpha_1\alpha_2 \\\quad+ [i(D_{11}-E_{11})+(D_{12}+E_{12})]^2 = 0.
\end{cases}
\]
Substituting $ h(z)=\alpha_1z_1 + \alpha_2z_2 + B $ into \eqref{eq-4.59}, we obtain
\begin{align}\label{eq-4.67}
	\dfrac{\partial f}{\partial z_1} - \dfrac{\partial f}{\partial z_2}&= (D_{11} - E_{11})\cos (\alpha_1z_1 + \alpha_2z_2 + B) \\&\nonumber - (D_{12} + E_{12})\sin (\alpha_1z_1 + \alpha_2z_2 + B) +(T_1 - T_2).
\end{align}
The characteristic equations of \eqref{eq-4.67} are $ \frac{dz_1}{dt}= 1,\; \frac{d z_2}{dt} = -1 $ and we see that
\begin{align}\label{eq-4.68}
	\dfrac{df}{dt}=(D_{11} - E_{11})\cos (\alpha_1z_1 + \alpha_2z_2 + B) - (D_{12} + E_{12})\sin (\alpha_1z_1 + \alpha_2z_2 + B) +(T_1 - T_2).
\end{align}
Using the initial conditions: $ z_1 = 0 $, $ z_2 = s $ and $ f =f(0,s):=\psi_{0}(s) $ with a parameter $ s $. Therefore, it is easy to see that $ z_1 = t,\; z_2 = -t + s, $
and
\begin{align*}
	f(t,s) &=\int_{0}^{t}[(D_{11} - E_{11})\cos (\alpha_1z_1 + \alpha_2z_2 + B) - (D_{12} + E_{12})\sin (\alpha_1z_1 + \alpha_2z_2 + B)]dt\\&\quad +\int_{0}^{t}(T_1 - T_2)dt + \psi_{0}(s)\\&= \int_{0}^{t}((D_{11} - E_{11})\cos [(\alpha_1 - \alpha_2)t + \alpha_2 s + B] - (D_{12} + E_{12})\sin [(\alpha_1 - \alpha_2)t + \alpha_2 s + B])dt\\&\quad +\int_{0}^{t}(T_1 - T_2)dt + \psi_{0}(s)\\&=\dfrac{(D_{11} - E_{11})}{(\alpha_1 - \alpha_2)}\sin [(\alpha_1 - \alpha_2)t + \alpha_2 s + B] + \dfrac{(D_{12} + E_{12})}{(\alpha_1 - \alpha_2)}\cos [(\alpha_1 - \alpha_2)t + \alpha_2 s + B]\\&\quad + (T_1 - T_2)t + \psi_{1}(s)
\end{align*}
where $ \psi_{1}(s) $ is a transcendental entire function with finite order in $ s $ such that
\begin{align*}
	\psi_{1}(s)=\psi_{0}(s) - \dfrac{(D_{11} - E_{11})}{(\alpha_1 - \alpha_2)}\sin (\alpha_2 s + B) - \dfrac{(D_{12} + E_{12})}{(\alpha_1 - \alpha_2)}\cos (\alpha_2 s + B).
\end{align*}
Therefore, we have
\begin{align}\label{eq-4.69}
	f(z_1,z_2)&=\dfrac{(D_{11} - E_{11})}{(\alpha_1 - \alpha_2)}\sin (\alpha_1z_1 + \alpha_2z_2 + B) + \dfrac{(D_{12} + E_{12})}{(\alpha_1 - \alpha_2)}\cos (\alpha_1z_1 + \alpha_2z_2 + B)\\&\nonumber + (T_1 - T_2)z_1 + \psi_{1}(s).
\end{align}
which can be written as
\begin{align*}
	f(z_1,z_2)=\dfrac{(D_{11} - E_{11})}{(\alpha_1 - \alpha_2)}\sin (L(z) + B) + \dfrac{(D_{12} + E_{12})}{(\alpha_1 - \alpha_2)}\cos (L(z)+ B) + (T_1 - T_2)z_1 + \psi_{1}(s).
\end{align*}
Substituting \eqref{eq-4.69} into \eqref{eq-4.57}, we can obtain that $ \psi_{1}(s) $ satisfies 
\begin{align*}
	&\psi^{\prime}_{1}(s) + \psi_{1}(s+s_0)\\ &= \left(\dfrac{E_{11}\alpha_1 - D_{11}\alpha_2}{\alpha_1-\alpha_2}\right)\cos (L(z)+B) + \left(\dfrac{E_{12}\alpha_1 - D_{12}\alpha_2}{\alpha_1-\alpha_2}\right)\sin (L(z)+B)\\& -\left(\dfrac{D_{11}-E_{11}}{\alpha_1-\alpha_2}\right)\sin(L(z)+L(c)+B) - \left(\dfrac{D_{12}-E_{12}}{\alpha_1-\alpha_2}\right)\cos(L(z)+L(c)+B)\\& - (T_1-T_2)(z_1+c) + T_2.
\end{align*}
This completes the proof of the Theorem \ref{th-3.4}.
\end{proof}	

\begin{proof}[\bf Proof of Theorem \ref{th-3.5}]
Suppose that $ f $ is a transcendental entire solution of \eqref{eq-3.5} with finite order. The equation \eqref{eq-3.5} can be written as
\begin{align}\label{eq-4.70}
	\left(\dfrac{DA^{\pm}}{-\Delta}\right)u^2 + \left(\dfrac{DB^{\mp}}{-\Delta}\right)v^2 = 1
\end{align}
Now we will discuss two cases below.
	
\noindent{\bf Case 1}: If $ \sqrt{\frac{DA^{\pm}}{-\Delta}}u $ is a constant, then it follows from \eqref{eq-4.70} that $ \sqrt{\frac{DB^{\mp}}{-\Delta}}v $ is also constant. Denote
\begin{align*}
	\sqrt{\dfrac{DA^{\pm}}{-\Delta}}u = K_9 ,\;\;\;\; \sqrt{\dfrac{DB^{\mp}}{-\Delta}}v =K_{10}
\end{align*}
In view of \eqref{eq-4.70}, it is easy to see that $ K^2_9 + K^2_{10} = 1 $. So we have 
\begin{align}\label{eq-4.71}
	M_1(f)= f(z+c) + \dfrac{\partial f}{\partial z_1}
	= \dfrac{K_9\xi^{\pm}_1}{\sqrt{\frac{DA^{\pm}}{-\Delta}}} - \dfrac{K_{10}\eta^{\pm}_1}{\sqrt{\frac{DB^{\mp}}{-\Delta}}} + \dfrac{\alpha\gamma- b\beta}{ab- \alpha^2}:= A_5,\\ \label{eq-4.72}	
	M_3(f) = f(z+c) + \dfrac{\partial^2 f}{\partial z^2_1}= \dfrac{K_9\eta^{\pm}_1}{\sqrt{\frac{DA^{\pm}}{-\Delta}}}+ \dfrac{K_{10}\xi^{\pm}_1}{\sqrt{\frac{DB^{\mp}}{-\Delta}}} + \dfrac{\alpha\beta- a\gamma}{ab- \alpha^2}:= B_5.
\end{align}
A simple computation using \eqref{eq-4.71} and \eqref{eq-4.72} shows that
\begin{align}\label{eq-4.73}
	\dfrac{\partial f}{\partial z_1} - \dfrac{\partial^2 f}{\partial z^2_1} =A_5 - B_5.
\end{align}
Differentiating \eqref{eq-4.71} with respect to $ z_1 $ we obtain
\begin{align}\label{eq-4.74}
	\dfrac{\partial f(z+c)}{\partial z_1} + \dfrac{\partial^2 f}{\partial z^2_1}=0
\end{align}
Thus it follows from \eqref{eq-4.72} and \eqref{eq-4.74} that
\begin{align}\label{eq-4.75}
	f(z+c) - \dfrac{\partial f(z+c)}{\partial z_1} = B_5
\end{align}
Differentiating \eqref{eq-4.75} with respect to $ z_1 $, we obtain
\begin{align}\label{eq-4.76}
	\dfrac{\partial f(z+c)}{\partial z_1} - \dfrac{\partial^2 f(z+c)}{\partial z^2_1}= 0
\end{align}
In view of \eqref{eq-4.73} and \eqref{eq-4.76} it is easy to see that $ A_5 = B_5 $. Using $ K^2_9 + K^2_{10} = 1 $, $ A_5 = B_5 $ and from \eqref{eq-4.71}, \eqref{eq-4.72}, a simple computation shows that $ K_2 =\frac{-R_2R_3\pm R_1\sqrt{R^2_1 +R^2_2 -R^2_3}}{R^2_1 + R^2_2} $, where 
\begin{align*}
	R_1 =\dfrac{\xi^{\pm}_1-\eta^{\pm}_1}{\sqrt{\frac{DA^{\pm}}{-\Delta}}},\; R_2=\dfrac{\xi^{\pm}_1+ \eta^{\pm}_1}{\sqrt{\frac{DB^{\mp}}{-\Delta}}}\; \mbox{and}\; R_3=\dfrac{\beta(\alpha + b)- \gamma(\alpha + a)}{ab-\alpha^2}.
\end{align*}
In view of \eqref{eq-4.73}, using $ A_5 = B_5 $, we obtain
\begin{align}\label{eq-4.77}
	f(z_1,z_2)= e^{z_1+\psi_1(z_2)} + \psi_{2}(z_2).
\end{align}
where $ \psi_1(z_2) $, $ \psi_2(z_2) $ are two function in $ z_2 $. From \eqref{eq-4.71} and \eqref{eq-4.77}, it follows that $ e^{(z_1+c_1)+\psi_{1}(z_2+c_2)} + e^{z_1+\psi_{1}(z_2)} + \psi_{2}(z_2+c_2) =A_5 $
which implies that $ \psi_{2}(z_2) =A_5,\;\; \psi_{1}(z_2)=Az_2 + B,\;\;\mbox{where}\;\; A=\frac{2k\pi i\pm\pi i-c_1}{c_2} $
Therefore, from \eqref{eq-4.77} we obtain $ f(z_1,z_2)= A_5 + e^{(z_1 + Az_2 + B)}. $\\
\noindent{\bf Case 2}:  If $ \sqrt{\frac{DA^{\pm}}{-\Delta}}u $ is not a constant. By the similar arguments being used in the proof of Theorem \ref{th-3.1}, there exists an entire function $ h(z) $ such that  
\begin{align}\label{eq-4.78}
	M_1(f)= f(z+c) + \dfrac{\partial f}{\partial z_1}=D_{11}\cos h(z) - D_{12}\sin h(z) + T_1, \\ \label{eq-4.79} M_3(f) = f(z+c) + \dfrac{\partial^2 f}{\partial z^2_1}=E_{11}\cos h(z) + E_{12}\sin h(z) + T_2.
\end{align}
In view of \eqref{eq-4.78} and \eqref{eq-4.79},  we obtain 
\begin{align}\label{eq-4.80}
	\dfrac{\partial f}{\partial z_1} - \dfrac{\partial^2 f}{\partial z^2_1}= (D_{11} - E_{11})\cos h(z) - (D_{12} + E_{12})\sin h(z) +(T_1 - T_2).
\end{align}
Differentiating \eqref{eq-4.78} with respect to $ z_1 $, we obtain
\begin{align}\label{eq-4.81}
	\dfrac{\partial f(z+c)}{\partial z_1} + \dfrac{\partial^2 f}{\partial z^2_1}=- D_{11}\sin h(z)\dfrac{\partial h}{\partial z_1} - D_{12}\cos h(z)\dfrac{\partial h}{\partial z_1}.
\end{align}
Thus it follows from \eqref{eq-4.79} and \eqref{eq-4.81} that
\begin{align}\label{eq-4.82}
	f(z+c) - \dfrac{\partial f(z+c)}{\partial z_1} =\left(D_{12}\dfrac{\partial h}{\partial z_1} + E_{11}\right)\cos h(z) + \left(D_{11}\dfrac{\partial h}{\partial z_1} + E_{12}\right)\sin h(z) + T_2
\end{align}
Differentiating \eqref{eq-4.82} with respect to $ z_1 $ and using \eqref{eq-4.80}, we obtain
\begin{align*}
	&\left(D_{11}\dfrac{\partial^2 h}{\partial z^2_1}-D_{12}\left(\dfrac{\partial h}{\partial z_1}\right)^2 -E_{11}\dfrac{\partial h}{\partial z_1}\right)\sin h(z) + \left(D_{12}\dfrac{\partial^2 h}{\partial z^2_1}+D_{11}\left(\dfrac{\partial h}{\partial z_1}\right)^2 +E_{12}\dfrac{\partial h}{\partial z_1}\right)\cos h(z)\\&\nonumber= (D_{11}-E_{11})\cos h(z+c) - (D_{12} + E_{12})\sin h(z+c) +(T_1 - T_2).
\end{align*}
This can be expressed as the following
\begin{align}\label{eq-4.83}
	&\left(\dfrac{\{(D_{12}-iD_{11})+(E_{12}+iE_{11})\}e^{ih(z+c)} + 2i(T_2 - T_1)}{(iD_{11} + D_{12})-(iE_{11}-E_{12})}\right)e^{ih(z+c)}\\&\nonumber + \left(\dfrac{(D_{11}+iD_{12})\dfrac{\partial^2 h}{\partial z^2_1} + (iD_{11}-D_{12})\left(\dfrac{\partial h}{\partial z_1}\right)^2 -(E_{11}-iE_{12})\dfrac{\partial h}{\partial z_1}}{(iD_{11} + D_{12})-(iE_{11}-E_{12})}\right)e^{i[h(z+c)+h(z)]}\\&\nonumber + \left(\dfrac{(iD_{12}-D_{11})\dfrac{\partial h}{\partial z_2} + (iD_{11}+D_{12})\left(\dfrac{\partial h}{\partial z_1}\right)^2 +(iE_{12}+E_{11})\dfrac{\partial h}{\partial z_1}}{(iD_{11} + D_{12})-(iE_{11}-E_{12})}\right)e^{i[h(z+c)-h(z)]} = 1.
\end{align}
If $ \left(\frac{(D_{11}+iD_{12})\frac{\partial^2 h}{\partial z^2_1} + (iD_{11}-D_{12})\left(\frac{\partial h}{\partial z_1}\right)^2 -(E_{11}-iE_{12})\frac{\partial h}{\partial z_1}}{(iD_{11} + D_{12})-(iE_{11}-E_{12})}\right) = 0 $, then from \eqref{eq-4.83} it is easy to see that $ \left(\frac{(iD_{12}-D_{11})\frac{\partial h}{\partial z_2} + (iD_{11}+D_{12})\left(\frac{\partial h}{\partial z_1}\right)^2 +(iE_{12}+E_{11})\frac{\partial h}{\partial z_1}}{(iD_{11} + D_{12})-(iE_{11}-E_{12})}\right)\neq 0 $. Otherwise, it follows that $ \left(\frac{\{(D_{12}-iD_{11})+(E_{12}+iE_{11})\}e^{ih(z+c)} + 2i(T_2 - T_1)}{(iD_{11} + D_{12})-(iE_{11}-E_{12})}\right) = 1 $, thus, a simple computation shows that $ e^{ih(z+c)}=\frac{i(T_1 - T_2)\pm i\sqrt{(T_2-T_1)^2 + [(D_{12}-iD_{11})+(E_{12}+iE_{11})][(iE_{11}-E_{12})-(iD_{11}+D_{12})]}}{(D_{12}-iD_{11})+(E_{12}+iE_{11})} $, which implies $ h(z) $ is a constant. Thus, the equation \eqref{eq-4.83} reduces to 
\begin{align*}
&\left(\dfrac{\{(D_{12}-iD_{11})+(E_{12}+iE_{11})\}e^{ih(z+c)} + 2i(T_2 - T_1)}{(iD_{11} + D_{12})-(iE_{11}-E_{12})}\right)e^{ih(z+c)}\\&\nonumber + \left(\dfrac{(iD_{12}-D_{11})\dfrac{\partial h}{\partial z_2} + (iD_{11}+D_{12})\left(\dfrac{\partial h}{\partial z_1}\right)^2 +(iE_{12}+E_{11})\dfrac{\partial h}{\partial z_1}}{(iD_{11} + D_{12})-(iE_{11}-E_{12})}\right)e^{i[h(z+c)-h(z)]} = 1.
\end{align*}
By the similar argument being used previously, in view of the Second Fundamental Theorem of Nevanlinna for several complex variables, we obtain $ T(r,e^{2ih(z+c)})=o\left(T(r,e^h)\right) $ which is a contradiction. Hence, $ \frac{(D_{11}+iD_{12})\frac{\partial^2 h}{\partial z^2_1} + (iD_{11}-D_{12})\left(\frac{\partial h}{\partial z_1}\right)^2 -(E_{11}-iE_{12})\frac{\partial h}{\partial z_1}}{(iD_{11} + D_{12})-(iE_{11}-E_{12})}\neq 0 $. Similarly, we obtain  $ \frac{(iD_{12}-D_{11})\frac{\partial h}{\partial z_2} + (iD_{11}+D_{12})\left(\frac{\partial h}{\partial z_1}\right)^2 +(iE_{12}+E_{11})\frac{\partial h}{\partial z_1}}{(iD_{11} + D_{12})-(iE_{11}-E_{12})}\neq 0 $. Therefore, in view of the Lemma \ref{lem-4.3} using \eqref{eq-4.83}, we obtain
\begin{align}\label{eq-4.85}
	\left(\dfrac{(iD_{12}-D_{11})\dfrac{\partial h}{\partial z_2} + (iD_{11}+D_{12})\left(\dfrac{\partial h}{\partial z_1}\right)^2 +(iE_{12}+E_{11})\dfrac{\partial h}{\partial z_1}}{(iD_{11} + D_{12})-(iE_{11}-E_{12})}\right)e^{i[h(z+c)-h(z)]}\equiv 1.
\end{align}	
From \eqref{eq-4.83} and \eqref{eq-4.85} it follows that
\begin{align}\label{eq-4.86}
	\left(\dfrac{(D_{11}+iD_{12})\dfrac{\partial^2 h}{\partial z^2_1} + (iD_{11}-D_{12})\left(\dfrac{\partial h}{\partial z_1}\right)^2 -(E_{11}-iE_{12})\dfrac{\partial h}{\partial z_1}}{(iD_{11} + D_{12})-(iE_{11}-E_{12})}\right)e^{i[h(z)-h(z+c)]}\equiv 1.
\end{align}
Since $ h(z) $ is a polynomial, then equation \eqref{eq-4.85}(or \eqref{eq-4.86}) implies $ h(z+c)-h(z)=\eta $, where $ \eta $ is a constant in $ \mathbb{C} $. Therefore, it follows that $ h(z)=L(z) + H(z) + B $, where $ L(z)=\alpha_1z_1 + \alpha_2z_2 $, $ H(z):= H(s_1) $, $ H(s_1) $ is a polynomial in $ s_1=c_2z_1 - c_1z_2 $. Now, we show that $ H(z)\equiv 0 $. From \eqref{eq-4.63} and \eqref{eq-4.64}, it follows that
\begin{align}\label{eq-4.87}
	&(iD_{11}+ D_{12})\alpha^2_1 + (iE_{12}+E_{11})\alpha_1 +c_2[2\alpha_1(iD_{11}+ D_{12})+ (iE_{12}+E_{11})]H^{\prime} \\&\nonumber + c^2_2(iD_{11}+ D_{12})(H^{\prime})^2 = e^{-iL(c)}[(iD_{11} + D_{12})-(iE_{11}-E_{12})]
\end{align}
and
\begin{align}\label{eq-4.88}
	&(iD_{11}-D_{12})\alpha^2_1 - (E_{11}-iE_{12})\alpha_1 +c_2[2\alpha_1(iD_{11}-D_{12})- (E_{11}-iE_{12})]H^{\prime} \\&\nonumber + c^2_2(iD_{11}-D_{12})(H^{\prime})^2 = e^{iL(c)}[(iD_{11} + D_{12})-(iE_{11}-E_{12})].
\end{align}
Combining with the fact that $ (c_1,c_2)\neq (0,0) $, it follows that $ H^{\prime} $ is a constant, that is $ \deg_{s_1}H\leq 1 $. Hence, $ L(z) + H(z) + B $ is still linear form of $ \alpha_1z_1 + \alpha_2z_2 + B $, which means that $ H(z)\equiv 0 $. Therefore it is easy to see that $ h(z)=L(z) + B =\alpha_1z_1 + \alpha_2z_2 + B $. Thus, the equations \eqref{eq-4.87} and \eqref{eq-4.88} become
\begin{align}\label{eq-4.89}
	\dfrac{(iD_{11}+ D_{12})\alpha^2_1 + (iE_{12}+E_{11})\alpha_1}{(iD_{11} + D_{12})-(iE_{11}-E_{12})}e^{iL(c)} = 1,
\end{align}
and
\begin{align}\label{eq-4.90}
	\dfrac{(iD_{11}-D_{12})\alpha^2_1 - (E_{11}-iE_{12})\alpha_1}{(iD_{11} + D_{12})-(iE_{11}-E_{12})}e^{-iL(c)} = 1.
\end{align}
From \eqref{eq-4.89} and \eqref{eq-4.90}, a simple computation shows that
\begin{align*}
	e^{2iL(c)} = \dfrac{(iD_{11}-D_{12})\alpha_1 - (E_{11}-iE_{12})}{(iD_{11}+ D_{12})\alpha_1 + (iE_{12}+E_{11})}
\end{align*}
and
\begin{align*}
	(D^2_{11} + D^2_{12})\alpha^4_1 &+ 2(D_{11}E_{12} +D_{12}E_{11})\alpha^3_1 + (E^2_{11} + E^2_{12})\alpha^2_1 \\& + [(iD_{11} + D_{12})-(iE_{11}-E_{12})]=0.
\end{align*}
The equation \eqref{eq-4.82} can be written as 
\begin{align*}
	f(z+c) - \dfrac{\partial f(z+c)}{\partial z_1} &=\left(D_{12}\alpha_1 + E_{11}\right)\cos (\alpha_1z_1 + \alpha_2z_2 + B) \\&\nonumber + \left(D_{11}\alpha_1 + E_{12}\right)\sin (\alpha_1z_1 + \alpha_2z_2 + B) + T_2,
\end{align*}
that is 
\begin{align}\label{eq-4.91}
	\dfrac{\partial f}{\partial z_1} - f(z) &=-\left(D_{12}\alpha_1 + E_{11}\right)\cos [L(z)-L(c)+B] \\&\nonumber - \left(D_{11}\alpha_1 + E_{12}\right)\sin [L(z)-L(c)+B] - T_2, 
\end{align}
Therefore, a simple computation shows that
\begin{align}\label{eq-4.92}
	f(z_1,z_2)= \phi(z_2)e^{z_1} + R_{11}\sin [L(z)-L(c)+B] + R_{12}\cos [L(z)-L(c)+B],
\end{align}
where $ \phi(z_2) $ is an entire function in $ z_2 $ and $ 	R_{11}=\frac{(D_{11}\alpha_1 + E_{12}) - \alpha_1(D_{12}\alpha_1 + E_{11})}{\alpha^2_1 + 1},\;\; R_{12}= \frac{(D_{12}\alpha_1 + E_{11}) + \alpha_1(D_{11}\alpha_1 + E_{12})}{\alpha^2_1 + 1} $. Substituting \eqref{eq-4.92} into \eqref{eq-4.78}, it is easy to see that
\begin{align*}
	&\phi(z_2 + c_2)e^{c_1}e^{z_1} +(R_{11}+D_{12})\sin (L(z)+B) + (R_{12}-D_{11})\cos (L(z)+B) \\& + R_{11}\alpha_1\cos [L(z)-L(c)+B] -R_{12}\alpha_1\sin [L(z)-L(c)+B] = T_1 - \phi(z_2)e^{z_1}. 
\end{align*}
Comparing both side the coefficient of $ e^{z_1} $ we obtain $ \phi(z_2 + c_2)= -e^{-c_1}\phi(z_2), $
which implies that $ \phi(z_2)=(-e^{-c_1})^{\frac{z_2}{c_2}}\pi(z_2),\;\;\mbox{where}\;\; \pi(z_2+c_2)=\pi(z_2). $ Therefore, from \eqref{eq-4.92}, we see that the precise form of the solution is 
\begin{align*}
	f(z_1,z_2)=(-e^{-c_1})^{\frac{z_2}{c_2}}\pi(z_2)e^{z_1} + R_{11}\sin [L(z)-L(c)+B] + R_{12}\cos [L(z)-L(c)+B].
\end{align*}
This completes the proof of the Theorem \ref{th-3.5}.
\end{proof}	

\noindent\textbf{Acknowledgment:} The authors would like to thank the referee for their helpful suggestions and comments to improve the exposition of the paper.
\vspace{1.6mm}

\noindent\textbf{Compliance of Ethical Standards:}\\

\noindent\textbf{Conflict of interest.} The authors declare that there is no conflict  of interest regarding the publication of this paper.\vspace{1.5mm}

\noindent\textbf{Data availability statement.}  Data sharing is not applicable to this article as no datasets were generated or analyzed during the current study.

\end{document}